\documentclass[11pt,leqno]{article}
\usepackage[centertags]{amsmath}
\usepackage{amsfonts}
\usepackage{amsmath}
\usepackage{amssymb}
\usepackage{latexsym}
\usepackage{amsthm}
\usepackage{mathrsfs}
\usepackage{newlfont}
\usepackage[latin1]{inputenc}

\newcommand{\D}{\displaystyle}
\newcommand{\Pn}{\{P_n\}_{n\geq0}}

\numberwithin{equation}{section} 

\title{Matrix valued orthogonal polynomials related to SU$(N+1)$,
their algebras of differential operators and the corresponding
curves \footnote{The
  work of the first author is  partially supported by  D.G.E.S,
  ref. BFM2003-06335-C03-01, FQM-262 (\textit {Junta de
Andaluc\'{\i}a}), that of the second author is partially supported
by NSF grant DMS-0204682.}}
\author{F. Alberto
Gr\"unbaum$^{\dagger}$ and  Manuel D. de la Iglesia$^{\ddagger}$ \\
   \footnotesize $\dagger$ \footnotesize\ Department of Mathematics. University of
California,  Berkeley \\
    \footnotesize Berkeley, CA 94720  U.S.A. grunbaum@math.berkeley.edu \\
   $\ddagger$  \footnotesize
    \  Departamento de An\'{a}lisis Matem\'{a}tico.
   Universidad de Sevilla \\
   \footnotesize Apdo (P. O. BOX) 1160. 41080 Sevilla. Spain.
mdi29@us.es
\\
\ \ }
\date{}

\begin{document}

\maketitle

\begin{abstract}
We give two examples of algebras of differential operators
associated to families of matrix valued orthogonal polynomials
arising from representations of SU$(N+1)$. The first one gives a
commutative algebra and the second one a non--commutative one.
\end{abstract}

\section{Introduction and outline of the paper}

The study of certain representations of U$(N)$, considered in
\cite{GPT1} for $N=2$ and fully explored in \cite{PT}, leads to a
second order differential operator $D$ with matrix valued
coefficients and a family of orthogonal matrix valued polynomials
$\Pn$  that are common eigenfunctions of this ordinary differential
operator with a {\em matrix valued} eigenvalue $\Lambda_n$,
$$
DP_n^*=P_n^*\Lambda_n.
$$

Different representations of U$(N)$ give rise to different
differential operators, each one with its corresponding families of
matrix valued orthogonal polynomials. Thus one obtains large sets of
examples of the type contemplated in \cite{D1}. A detailed look at
two such examples arising from the representation theory of U$(N)$
is given in \cite{GPT6}. The theory of matrix valued orthogonal
polynomials, without any consideration of differential equations,
goes back to \cite{K1,K2}. The question raised in \cite{D1} is a
matrix version of the question raised and settled in the scalar case
in \cite{B}, see also \cite{R}. By considering this non--commutative
setup one obtains a situation which is \textbf{much richer} than the
scalar one. For a sample of applications of matrix valued orthogonal
polynomials, see \cite{SV}.

Starting in  \cite{CG2} one considers {\em for a fixed family of
matrix valued orthogonal polynomials} the algebra of all such
differential operators $D$. Each differential operator in the
algebra will bring in its own sequence of matrix valued eigenvalues,
which can be denoted by $\Lambda_n(D)$. This algebra of matrix
valued differential operators and the algebra of matrix valued
eigenvalues are isomorphic, as we will show in subsection
\ref{sec22}.

In \cite{CG2} one finds a discussion of this general question in a
few instances where the family of matrix valued orthogonal
polynomials $\Pn$ is not related to any underlying group. These
examples were obtained previously in \cite{DG1} and \cite{G}.

In section \ref{sec1} we introduce two differential operators
obtained by conjugating another pair of "original" differential
operators arising from representation theory in \cite{GPT6}. Section
3 is devoted to the first one and section 4 to the second one.

In subsection \ref{sec21} we give another conjugation of one of the
original operators and compute a sequence of orthogonal polynomials
which are common eigenfunctions of it. Those polynomials are then
related to those obtained using the new tools introduced in
\cite{T2} in terms of the matrix version of Gauss' equation.

In subsection \ref{sec22} we return to the conjugation in section
\ref{sec1} and compute its corresponding family of orthogonal
polynomials and consider an algebra of differential operators
attached to this sequence. This family of orthogonal polynomials is
fully explored in the appendix. In subsections \ref{sec31} and
\ref{sec32}, this process is repeated for the second original
operator

Using these families detailed in the appendix, we present convincing
computational evidence that leads us to state the following: in the
first example the algebra of operators is commutative, while in the
second example this is not so. In each case, we can exhibit the
generators of the algebra and a basis for the space of operators of
a given order. Finally, in the first example we can write down
explicitly a polynomial relation among the \textbf{two} generators
of the algebra, while in the second example we give a collection of
some of the relations among the \textbf{five} generators of the
algebra.

\section{From the Casimir operator to some examples of ordinary
differential operators with matrix coefficients}\label{sec1}

In this section we give an extremely  sketchy indication of the way
in which one goes from the Casimir operator, namely a differential
operator acting on matrix valued functions on the group SU$(N+1)$,
to certain ordinary differential operators acting on matrix valued
functions defined in the interval $[0,1]$. There are many details
missing in the presentation below and the reader is advised to
consult either \cite {GPT6} or even better \cite{PT, GPT1} for a
full account.

The main goal here is to give the explicit differential operators
that will play a crucial role in the rest of the paper.

In this paper, following the lines of \cite{GPT6}, we use two
different skew hermitian matrix valued forms, namely
$$
(P,Q)=\int_{\mathbb{R}} P(t)\,W(t)\,Q^*(t)\,dt
$$
and
$$
\langle P,Q\rangle=(P^*,Q^*)^*.
$$
Here, $W(t)$ is a positive definite matrix valued weight.

Two matrix valued functions are called \emph{orthogonal} if
$(P,Q)=\Theta$, the null matrix of the appropriate dimension. On the
other hand, and for reasons explained in \cite{GPT6}, a differential
operator $D$ is called \emph{symmetric} if
$$
\langle D P,Q\rangle=\langle P,D Q\rangle,
$$
for all matrix valued polynomials $P$ and $Q$.

In the first example, considering representations $\pi$ of GL($N
,\mathbb{C}$), that correspond to \emph{one step} $N$--tuples of the
form
$$
\pi=(\underbrace{m+2,\ldots,m+2}_k,\underbrace{m,\ldots,m}_{N-k}),\quad
1\leq k\leq N-1
$$
and making the changes $N=\beta+1$ and $m=\alpha$, we have one
family of examples that depends on three parameters $\alpha,\beta,k$
where the differential operator is given by:

\begin{align*}
D = t(1-t) \frac {d^2}{dt^2} &+  \left[ \begin{pmatrix}
\alpha+3 &0 &0 \\
0 &\alpha+2 &0 \\
0 &0 &\alpha+1
\end{pmatrix} - t \begin{pmatrix}
\alpha+\beta+4 &0 &0 \\
0 &\alpha+\beta+3&0 \\
0 &0 &\alpha+\beta+2
\end{pmatrix} \right] \frac {d}{dt}   \\
&+ \frac {1}{1-t} \begin{pmatrix}
-2(\beta-k+1) &2(\beta-k+1) &0 \\
0 &-(\beta-k+2) &\beta-k+2 \\
0 &0 &0
\end{pmatrix} \displaybreak[0]\\
&\quad  + \frac {t}{1-t} \begin{pmatrix}
0 &0 &0 \\
k+1 &-(k+1) &0\\
0 & 2k&-2k
\end{pmatrix}.
\end{align*}
This operator is symmetric with respect to the matrix weight
function

$$W(t)=t^{\alpha}(1-t)^{\beta}\begin{pmatrix}
  w_1t^2 & 0 & 0 \\
  0 & w_2t & 0 \\
  0 & 0 & w_3
\end{pmatrix},\; \alpha,\beta>-1,$$
where
\begin{equation*}
   \begin{split}
      w_1&=\prod_{j=0}^{k-2}\frac{(\alpha-j-1)(\alpha-j)}{(k-j-1)(k-j)}\\
      w_2&=k(\alpha-k)\prod_{j=0}^{k-2}\frac{(\alpha-j-1)(\alpha-j)}{(k-j)(k-j+1)}\\
      w_3&=\prod_{j=0}^{k-1}\frac{(\alpha-j-1)(\alpha-j)}{(k-j)(k-j+1)}.
    \end{split}
\end{equation*}
The example above arises naturally in the context of representation
theory with $\beta$ and $k$ natural numbers. It is however possible
to consider both $D$ and $W(t)$ under the conditions $\beta>-1$ and
$1\leq k\leq\beta$.

The name {\em one-step} for the $(\beta+1)$--tuples or partitions
discussed above is very natural if one looks at the corresponding
Young diagrams, see \cite{VK}. The same is true in the second type
of example, discussed later in this section.

This operator is not written in {\em hypergeometric form} yet (the
coefficients are not matrix polynomials of degree less than or equal
to the corresponding order of differentiation). Just as in
\cite{GPT4} we proceed to do an appropriate conjugation by a certain
matrix valued function, namely
$$\Psi^*(t)=\begin{pmatrix}
  1 & 0 & 0 \\
  1 & 1 & 0\\
  1 & 2 & 1
\end{pmatrix}\begin{pmatrix}
  1 & 0 & 0 \\
  0 & 1-t & 0\\
  0 & 0 & (1-t)^2
\end{pmatrix} $$
such that we get$$\widetilde{W}(t)=\Psi(t)W(t)\Psi^*(t)$$ for the
new weight function. The new differential operator becomes $\tilde
DF=(\Psi^*)^{-1} D(\Psi^* F)$, where
\begin{equation*}
   \widetilde{D}= \widetilde{A}_2(t)\frac{d^2}{dt^2}+
\widetilde{A}_1(t)\frac{d}{dt}+ \widetilde{A}_0(t),
\end{equation*}
with $\widetilde{A}_2$, $ \widetilde{A}_1$, $ \widetilde{A}_0$ given
by:
\begin{eqnarray*}
\widetilde{A}_2(t) &=  &t(1-t),  \\
\widetilde{A}_1(t) &=  &\begin{pmatrix}
\alpha+3 &0 &0 \\
-1 & \alpha+2 &0 \\
0 & -2 &\alpha+1
\end{pmatrix} - t \begin{pmatrix}
\alpha+\beta+4 & 0&0 \\
0 &\alpha+\beta+5 &0 \\
0 &0 &\alpha+\beta+6
\end{pmatrix},  \\
\widetilde{A}_0(t) &= &\begin{pmatrix}
0 &2(\beta-k+1) &0 \\
0 &-(\alpha+\beta-k+2) &\beta-k+2 \\
0 &0 &-2(\alpha+\beta-k+3)
\end{pmatrix},
\end{eqnarray*}
written now in hypergeometric form. This differential operator is
symmetric with respect to the new weight function $\widetilde{W}(t)$
and will be the starting point of the discussion in subsection
\ref{sec22}.

As a remark connecting the matrix weight given above with some
considerations in \cite{DG1,PT1}, note that $\widetilde{W}(t)$
admits the factorization
$$
\widetilde{W}(t)=\frac{\rho(t)}{\rho(1/2)}T(t)\widetilde{W}(1/2)T^*(t)
$$
with $T(1/2)=I$ and $\rho(t)=t^{\alpha}(1-t)^{\beta}$. The matrix
$T(t)$ solves the equation
$$
T'(t)=\bigg(\frac{A}{t}+\frac{B}{1-t}\bigg)T(t)
$$
with
$$
A=\begin{pmatrix}
  1 & -\frac{1}{2} & 0 \\
  0 & \frac{1}{2} & -1 \\
  0 & 0 & 0
\end{pmatrix},
$$
$$
B=\begin{pmatrix}
  0 & -\frac{1}{2} & 0 \\
  0 & -1 & -1 \\
  0 & 0 & -2
\end{pmatrix}.
$$

To obtain  the second example we consider representations $\pi$ of
GL($N ,\mathbb{C}$) that correspond to \emph{two steps} $N$--tuples
of the form
$$
\pi=(\underbrace{m+2,\ldots,m+2}_{k_1},\underbrace{m+1,\ldots,m+1}_{k_2-k_1},\underbrace{m,\ldots,m}_{N-k_2}),\quad
1\leq k_1<k_2\leq N-1
$$
and replacing $N$ and $m$ as above, we have one family of examples
that depends on four parameters $\alpha,\beta,k_1,k_2$ where the
differential operator is given by:

\begin{align*}
D&=  t(1-t) \frac {d^2}{dt^2} + \left[
\begin{pmatrix}
\alpha+3 &0 &0 &0\\
0 &\alpha+2 &0 &0\\
0 &0 &\alpha+2 &0\\
0&0&0&\alpha+1
\end{pmatrix} \right. \\
   &  \qquad \qquad \qquad \quad\qquad \left.- t \begin{pmatrix}
\alpha+\beta+4 &0 &0 &0\\
0 & \alpha+\beta+3 &0 &0\\
0 & 0 & \alpha+\beta+3 &0\\
0 & 0 & 0 &\alpha+\beta+2
\end{pmatrix} \right] \frac {d}{dt} \displaybreak[0] \\
&\quad + \frac {1}{1-t} \begin{pmatrix} k_1+k_2-2(\beta+1)
&\tfrac{(k_2-k_1+2)(\beta-k_2+1)}{k_2-k_1+1} &
\tfrac{(k_2-k_1)(\beta-k_1+2)}{k_2-k_1+1}&0  \\
0 &-(\beta-k_1+2) &0 &\beta-k_1+2 \\
0&0& -(\beta-k_2+1)& \beta-k_2+1\\
   0 &0 &0 &0
\end{pmatrix}\displaybreak[0] \\
&\quad +  \frac {t}{1-t} \begin{pmatrix}
0 &0 &0 &0 \\
k_2+1 &-(k_2+1) &0 &0\\
k_1 & 0& -k_1&0 \\
0& \tfrac{k_1(k_2-k_1+2)}{k_2-k_1+1} &
\tfrac{(k_2-k_1)(k_2+1)}{k_2-k_1+1} & -(k_1+k_2)
\end{pmatrix}.
\end{align*}
This operator is symmetric with respect to the matrix weight
function

$$W(t)=t^{\alpha}(1-t)^{\beta}\begin{pmatrix}
  w_1t^2 & 0 & 0 &0\\
  0 & w_2t & 0&0 \\
  0 & 0 & w_3t&0\\
  0 & 0 & 0&w_4
\end{pmatrix},\; \alpha,\beta>-1,$$
where
\begin{equation*}
   \begin{split}
      w_1&=\frac{k_2-k_1+1}{k_2}\begin{pmatrix}
        \beta \\
        k_2-1 \\
      \end{pmatrix}\begin{pmatrix}
        \beta+1\\
        k_1-1 \\
      \end{pmatrix}\\
      w_2&=\frac{k_2-k_1+2}{k_2+1}\begin{pmatrix}
        \beta \\
        k_2 \\
      \end{pmatrix}\begin{pmatrix}
        \beta+1\\
        k_1-1 \\
      \end{pmatrix}\\
      w_3&=\frac{k_2-k_1}{k_2}\begin{pmatrix}
        \beta \\
        k_2-1 \\
      \end{pmatrix}\begin{pmatrix}
        \beta+1\\
        k_1 \\
      \end{pmatrix}\\
      w_4&=\frac{k_2-k_1+1}{k_2+1}\begin{pmatrix}
        \beta \\
        k_2 \\
      \end{pmatrix}\begin{pmatrix}
        \beta+1\\
        k_1 \\
      \end{pmatrix}
    \end{split}
\end{equation*}
Once again, both $D$ and $W(t)$ can be considered for $\beta>-1$ and
$1\leq k_1<k_2\leq \beta$.

As in the previous example, this operator is not written in {\em
hypergeometric form} yet. A possible conjugation in this case is
given by the matrix valued function, see \cite{GPT6}
$$\Psi^*(t)=\begin{pmatrix}
  1 & 0 & 0 &0\\
  1 & 1 & 0&0\\
  1 & 0 & 1 & 0\\
  1 & \frac{k_2-k_1+2}{k_2-k_1+1} & \frac{k_2-k_1}{k_2-k_1+1} & 1
\end{pmatrix}\begin{pmatrix}
  1 & 0 & 0 &0 \\
  0 & 1-t & 0 &0\\
  0 & 0 & 1-t &0\\
  0 & 0 & 0 & (1-t)^2
\end{pmatrix}, $$
and we get$$\widetilde{W}(t)=\Psi(t)W(t)\Psi^*(t)$$ for the new
weight function. The new differential operator becomes $\tilde
DF=(\Psi^*)^{-1} D(\Psi^* F)$, where
\begin{equation*}
   \widetilde{D}= \widetilde{A}_2(t)\frac{d^2}{dt^2}+
\widetilde{A}_1(t)\frac{d}{dt}+ \widetilde{A}_0(t),
\end{equation*}
with $\widetilde{A}_2$, $ \widetilde{A}_1$, $ \widetilde{A}_0$ given
by:
\begin{eqnarray*}
\widetilde{A}_2(t) &= &t(1-t) \\
\widetilde{A}_1(t) &=  &\begin{pmatrix}
\alpha+3 &0 &0 &0\\
-1 & \alpha+2 &0 &0 \\
-1 &0 & \alpha+2 &0\\
0 & -\tfrac{k_2-k_1+2}{k_2-k_1+1} & -\tfrac{k_2-k_1}{k_2-k_1+1}
&\alpha+1
\end{pmatrix} \\
& &\quad - t \begin{pmatrix}
\alpha+\beta+4 & 0&0 &0 \\
0 &\alpha+\beta+5 &0 &0\\
0& 0& \alpha+\beta+5&0\\
0&0 &0 &\alpha+\beta+6
\end{pmatrix}  \\
\widetilde{A}_0(t) &= &\begin{pmatrix} 0 &
\tfrac{(k_2-k_1+2)(\beta-k_2+1)}{k_2-k_1+1} &
\tfrac{(k_2-k_1)(\alpha-k_1+2)}{k_2-k_1+1}&0 \\
0 &-(\alpha+\beta+2)+k_2 &0 &\beta-k_1+2 \\
0&0& -(\alpha+\beta+3)+k_1 & \beta-k_2+1\\
0&0 &0 &-2(\alpha+\beta+3)+k_1+k_2
\end{pmatrix}.
\end{eqnarray*}
written now in hypergeometric form. This differential operator is
symmetric with respect to the new weight function $\widetilde{W}(t)$
and will be the starting point of the discussion in subsection
\ref{sec32}.

Again, as a remark connecting the matrix weight given above with
some considerations in \cite{DG1,GPT6}, note that $\widetilde{W}(t)$
admits the factorization
$$
\widetilde{W}(t)=\frac{\rho(t)}{\rho(1/2)}T(t)\widetilde{W}(1/2)T^*(t)
$$
with $T(1/2)=I$ and $\rho(t)=t^{\alpha}(1-t)^{\beta}$. The matrix
$T(t)$ solves the equation
$$
T'(t)=\bigg(\frac{A}{t}+\frac{B}{1-t}\bigg)T(t)
$$
with
$$
A=\begin{pmatrix}
  1 & -\frac{1}{2} & -\frac{1}{2} & 0 \\
  0 & \frac{1}{2} & 0 & \frac{k_1-k_2-2}{2(k_2-k_1+1)} \\
  0 & 0 & \frac{1}{2} & \frac{k_1-k_2}{2(k_2-k_1+1)} \\
  0 & 0 & 0 & 0
\end{pmatrix},
$$
$$
B=\begin{pmatrix}
  0 & -\frac{1}{2} & -\frac{1}{2} & 0 \\
  0 & -1 & 0 & \frac{k_1-k_2-2}{2(k_2-k_1+1)} \\
  0 & 0 & -1 & \frac{k_1-k_2}{2(k_2-k_1+1)} \\
  0 & 0 & 0 & -2
\end{pmatrix}.
$$

The point of these conjugations is to insure that the new
differential operator should have matrix valued coefficients that
are polynomials in $t$ of degrees not higher than the corresponding
order of differentiation. This new form of the differential operator
is referred to as {\bf a hypergeometric form}.

Given a differential operator there may be more than one way of
conjugating it into a hypergeometric form. In the next section we
will introduce another operator, denoted by $D$, which is obtained
from the one we called $D$ at the beginning of this section, by
means of a conjugation that is different from the one used to
produce $\widetilde{D}$ out of $D$. We will retain the symbol $W(t)$
to denote the new weight matrix obtained as in this section by
conjugation of the original $W(t)$. A similar  step will be taken in
section \ref{sec31}.

\section{The one step example}
\subsection{Generating the polynomial eigenfunctions of the differential
operator}\label{sec21}

We retain the symbol $D$ for the {\bf new} differential operator
given in {\em hypergeometric form} for the \emph{one step} example
of the previous section. It is now given by the expression:

$$
D=t(1-t)\D\frac{d^2}{dt^2}+(X-tU)\D\frac{d}{dt}+V,
$$
where
$$
V= \begin{pmatrix}
0 & 0 & 0 \\
0 & -(\alpha+\beta-k+2) & 0\\
0 & 0 & -2(\alpha+\beta-k+3)
\end{pmatrix},
$$
$$
U= \begin{pmatrix}
\alpha+\beta+4 & -1 & \frac{2}{\alpha+\beta-k+2}\\
0 & \alpha+\beta+5 & -\frac{2(\alpha+\beta-k+3)}{\alpha+\beta-k+2}\\
0 & 0 & \alpha+\beta+6
\end{pmatrix}\quad\mbox{and}
$$
$$
X= \begin{pmatrix}
\frac{(\alpha+1)(\alpha+\beta-k+4)}{\alpha+\beta-k+2} &
\frac{(\alpha+1)(\alpha+\beta-k+4)}{(\alpha+\beta-k+2)(\alpha+\beta
-k+3)} & 0\\
\frac{2(\beta-k+1)}{\alpha+\beta-k+2} &
\frac{C_1(\alpha,\beta,k)}{(\alpha+\beta-k+2)(\alpha+\beta-k
+4)} & \frac{2(\alpha+2)}{\alpha+\beta-k+4}\\
0 &
\frac{(\beta-k+2)(\alpha+\beta-k+2)}{(\alpha+\beta-k+3)(\alpha+\beta-k+4)}
& \frac{C_2(\alpha,\beta,k)}{\alpha+\beta-k+4}
\end{pmatrix},
$$
where \begin{eqnarray*}
\lefteqn{C_1(\alpha,\beta,k)=\alpha\beta^2+10\alpha\beta+18\alpha+2\beta^2+14\beta+16+2\alpha^2\beta+8\alpha^2}\\
\noalign{\smallskip} &
&\qquad\quad-2k\alpha\beta-10k\alpha-4k\beta-14k+\alpha^3-2k\alpha^2+k^2\alpha+2k^2
\end{eqnarray*} and
$$C_2(\alpha,\beta,k)=\alpha\beta+5\alpha+3\beta+8+\alpha^2-k\alpha-3k.$$

The conjugation that was used above to produce $D$ was chosen so
that $V$ turns out to be diagonal. The coefficients are not so
pleasant as the ones we will use in the next subsection, but it will
be easier to find a relation between some matrix valued orthogonal
polynomials whose adjoints are eigenfunctions of this new $D$ and
the \emph{matrix hypergeometric function}, introduced in \cite{T2}.

The operator  $D$ above and the new weight $W(t)$ are such that the
(unique) sequence of monic matrix valued polynomials
$\{Q_n\}_{n\geq0}$ with respect to $W(t)$ satisfies
$$
  D Q_n^*\equiv t(1-t)\D\frac{d^2}{dt^2}Q_n^*+(X-tU)\D\frac{d}{dt}Q_n^*+VQ_n^*=Q_n^*\Gamma_n,
$$
where $\Gamma_n=-n^2+n(I-U)+V$.

Now, if $\Pn$ is a sequence given by
$$
Q_n(t)=S_n^{-1}P_n(t),\quad \mbox{det}\; S_n\neq 0,
$$
we have that
$$
D P_n^*(t)=P_n^*(t)(S_n^*)^{-1}\Gamma_n S_n^*.
$$
We now make the genericity assumption that the eigenvalues of
$\Gamma_n$ are different for all $n$. Then by choosing $S_n^*$ as
the matrix whose columns are the eigenfunctions of $\Gamma_n$
(unique up to scaling) we have that $P_n^*(t)$ satisfy

\begin{equation}\label{eqdiff}
  D P_n^*\equiv t(1-t)\D\frac{d^2}{dt^2}P_n^*+(X-tU)\D\frac{d}{dt}P_n^*+VP_n^*=P_n^*\Lambda_n,
\end{equation}
with
$$
\Lambda_n= \begin{pmatrix}
t_1 & 0 & 0 \\
0 & t_2 & 0 \\
0 & 0 & t_3
\end{pmatrix}=(S_n^*)^{-1}\Gamma_n S_n^*
$$
and the values of $t_i$ are given in (\ref{diag}). If we put
$P_n^*(t)=\sum_{j=0}^n A_j^n t^j$, the equations satisfied by
$A_j^n$, $j=n,n-1,\cdots,0$, are:

\begin{equation}\label{coeff}
   \begin{split}
      \Gamma_nA_n^n&=A_n^n\Lambda_n\\
      \Gamma_jA_j^n-A_j^n\Lambda_n&=-(j+1)(X+j)A_{j+1}^n,\quad
      j=n-1,\cdots,0.
    \end{split}
\end{equation}
All these equations are of Sylvester´s type, see \cite{Gr}. Under
the genericity assumption made earlier, the first equation has a
solution $A_n^n$ that is unique up to the choice of three scalars.
Each one of the equations that follow has a unique solution $A_j^n$
if we assume an extra genericity assumption, namely for each $n$ the
spectrum of each $\Gamma_j$, $j=0,\cdots, n-1$ is disjoint from the
spectrum of $\Lambda_n$. For a careful treatment of the non generic
case, the reader can consult \cite{PR}.

If $\Pn$ is a family of matrix valued polynomials satisfying
(\ref{eqdiff}) and $D$ is symmetric it follows that
$$
\langle D P_n^*, P_m^*\rangle=\langle P_n^*\Lambda_n,
P_m^*\rangle=(P_m,P_n)\Lambda_n
$$
and
$$
\langle P_n^*,D P_m^*\rangle=\langle P_n^*,
P_m^*\Lambda_m\rangle=\Lambda_m^*(P_m,P_n).
$$

Now, under the genericity assumptions made earlier, the spectrum of
$\Lambda_n$ and that of $\Lambda_m^*$ are disjoint if $n\neq m$. The
classical uniqueness result for Sylvester´s equations gives that
$P_n$ and $P_m$ are orthogonal if $n\neq m$.

It is now our purpose to obtain a relation between the matrix valued
orthogonal polynomials $\Pn$ introduced above and the matrix
hypergeometric function. We use the tools in \cite{T2}. For a
warm--up, the reader can consult \cite{GPT5}, where these tools are
used in the same fashion as below. The main idea is to replace
elements in $M(3,\mathbb{C})$ by vectors in $\mathbb{C}^9$. We will
denote this map by vec. It will be important to replace right and
left multiplication in $M(3,\mathbb{C})$ by linear maps in
$\mathbb{C}^9$.

This allows us to rewrite the differential equation above as the
following equivalent differential equation:

$$
t(1-t)\D\frac{d^2}{dt^2}\mbox{vec}(P_n^*)+(C-t\tilde{U})\D\frac{d}{dt}\mbox{vec}(P_n^*)-\tilde{T}\mbox{vec}(P_n^*)=\Theta
$$
where $C$ and $\widetilde{U}$ are the $9\times9$ matrices obtained
by the rules, see \cite{HJ}:
$$
C=X\otimes I,\quad \widetilde{U}=U\otimes I,\quad\mbox{and}\quad
\widetilde{T}= V\otimes I -I\otimes \Lambda_n^*.
$$
One gets that $\widetilde{T}$ is given as follows:
$$
\tilde{T}=\mbox{diag}\{t_1,t_2,t_3,t_4,t_5,t_6,t_7,t_8,t_9\}
$$
where

\begin{equation}\label{diag}
   \begin{split}
      t_1&=-n^2-n(\alpha+\beta+3)\\
      t_2&=-n^{2} -n(\alpha+\beta+4)-(\alpha+\beta-k+2)\\
      t_3&=-n^{2} -n(\alpha+\beta+5)-2(\alpha+\beta-k+3)\\
      t_4&=-n^2-n(\alpha+\beta+3)+\alpha+\beta-k+2\\
      t_5&=-n^{2} -n(\alpha+\beta+4)\\
      t_6&=-n^{2} -n(\alpha+\beta+5)-(\alpha+\beta-k+4)\\
      t_7&=-n^2-n(\alpha+\beta+3)+2(\alpha+\beta-k+3)\\
      t_8&=-n^{2} -n(\alpha+\beta+4)+\alpha+\beta-k+4\\
      t_9&=-n^{2} -n(\alpha+\beta+5)
    \end{split}
\end{equation}
Continuing with the strategy in \cite{T2} we need to find matrices
$A$ and $B$ such that
$$\tilde{U}=I+A+B\quad \mbox{and}\quad
\tilde{T}=AB.$$ The fact that $\tilde{T}$ is diagonal will make this
easier than it would be otherwise. This was the reason for choosing
the conjugation that made $V$ diagonal. The introduction of $A$ and
$B$ allows us to rewrite our new equation in a form that is very
much like the {\bf classical hypergeometric equation of Euler and
Gauss}, see (\ref{eqd3}).

The factorization in the last equation is not unique, and we look
for $A$ in the form

$$
A=\begin{pmatrix}
  A_{11} & A_{12} & A_{13}  \\
  \Theta & A_{22} & A_{23} \\
  \Theta & \Theta & A_{33}
\end{pmatrix}
$$
where each $A_{ij}$ is a $3\times 3$ diagonal matrix. Once we have
found $A$, the matrix $B$ will be given by the expression
$$
B=\tilde{U}-A-I.
$$
With this form for $B$ the factorization above gives a number of
conditions on $A$.

We put
$$ A_{ii}=\begin{pmatrix}
  \phi_i & 0 & 0  \\
  0 & \varphi_i & 0  \\
  0 & 0 & \psi_i
\end{pmatrix}
\;\; i=1,2,3.
$$
with $\phi_i,\varphi_i,\psi_i\;i=1,2,3$ to be determined later.

The elements of the diagonal of $A_{12}$ are of the form:

$$
\D\frac{\gamma_1}{-\gamma_1-\gamma_2+\alpha+\beta+4},
$$
with $\gamma=\phi,\varphi,\psi$ in each entry.

If we denote by

$$
\omega_{13}=\D\frac{2}{\alpha+\beta-k+2}
$$
$$
\omega_{23}=-\D\frac{2(\alpha+\beta-k+3)}{\alpha+\beta-k+2}
$$
the elements of the last column of $U$, the elements of the diagonal
of $A_{23}$ turn out to be of the form:

$$
\D\frac{-\omega_{23}\gamma_2}{-\gamma_2-\gamma_3+\alpha+\beta+5},
$$
with $\gamma=\phi,\varphi,\psi$ in each entry.

The elements of the diagonal of $A_{13}$ are a bit more complicated,
they are given by

\begin{eqnarray*}
-\frac{\omega_{13}\gamma_1}{-\gamma_1-\gamma_3+\alpha+\beta+5}-\frac{\omega_{23}\gamma_1}
{(-\gamma_1-\gamma_3+\alpha+\beta+5)(-\gamma_1-\gamma_2+\alpha+\beta+4)}
\bigskip\\
-\frac{\omega_{23}\gamma_1\gamma_2}{(-\gamma_1-\gamma_3+\alpha+\beta+5)(-\gamma_1-\gamma_2+\alpha+\beta+4)
(-\gamma_2-\gamma_3+\alpha+\beta+5)},
\end{eqnarray*}

with $\gamma=\phi,\varphi,\psi$ in each entry.

The parameters $\phi_i,\varphi_i,\psi_i$, $i=1,2,3$ are subject to
the following conditions, resulting from the factorization above:

$$
\left\{\begin{array}{l}
\phi_1=-n\quad \mbox{or}\quad \phi_1=n+\alpha+\beta+3\\
\varphi_1^2-(\alpha+\beta+3)\varphi_1+t_2=0\\
\psi_1^2-(\alpha+\beta+3)\psi_1+t_3=0
\end{array}\right.
$$
$$
\left\{\begin{array}{l}
\phi_2^2-(\alpha+\beta+4)\phi_2+t_4=0\\
\varphi_2=-n\quad \mbox{or}\quad \varphi_2=n+\alpha+\beta+4\\
\psi_2^2-(\alpha+\beta+4)\psi_2+t_6=0
\end{array}\right.
$$
$$
\left\{\begin{array}{l}
\phi_3^2-(\alpha+\beta+5)\phi_3+t_7=0\\
\varphi_3^2-(\alpha+\beta+5)\varphi_3+t_8=0\\
\psi_3=-n\quad \mbox{or}\quad \psi_3=n+\alpha+\beta+5.
\end{array}\right.
$$

We have reached the main point of \cite{T2}, namely, the vector
$\mbox{vec}(P_n^*)$ satisfies the matrix hypergeometric equation:

\begin{equation}\label{eqd3}
t(1-t)\D\frac{d^2}{dt^2}\mbox{vec}(P_n^*)+(C-t(I+A+B))\D\frac{d}{dt}\mbox{vec}(P_n^*)-AB\mbox{vec}(P_n^*)=\Theta.
\end{equation}
Here, by abuse of notation, $\Theta$ denotes the null vector.

The eigenvalues of $C$ are $\{\alpha+1,\alpha+2,\alpha+3\}$ with
multiplicity 3 each one, so if $\alpha>-1$ we meet the conditions
required to make sense of the following definition:
$$
(C,A,B)_{i+1}=(C+iI)^{-1}(A+iI)(B+iI)(C+(i-1)I)^{-1}(A+(i-1)I)(B+(i-1)I)\cdots
C^{-1}AB
$$
for all $i\geq 0$ and $(C,A,B)_0=I$.

Then we define the\emph{ matrix hypergeometric function} as:
$$
_2F_1(C,A,B;t)=\D\sum_{i\geq0}(C,A,B)_i\D\frac{t^i}{i!}.
$$

It is proved in \cite{T2} that $_2F_1(C,A,B;t)$ is analytic on
$|t|<1$, and the analytic solutions at $t=0$ of (\ref{eqd3}) are
given by
$$
\mbox{vec}\bigl(\Phi(t)\bigr)=
\,_2F_1(C,A,B;t)\mbox{vec}\bigl(\Phi(0)\bigr).
$$

The function $_2F_1(C,A,B;t)$ is not a polynomial function, as in
the classical case, but nevertheless we have that
\begin{equation}\label{vec}
\mbox{vec}\bigl(P_n^*(t)\bigr)=\,_2F_1(C,A,B;t)\mbox{vec}\bigl(P_n^*(0)\bigr)
\end{equation}
is a vector valued polynomial of degree $n$ in $t$.

In (\ref{vec}) we can give explicitly the value of $P_n^*(0)$:
following the strategy explained around (\ref{coeff}), we have
$P_n^*(0)=A_0^n$, and we get:

$$
P_n^*(0)=\begin{pmatrix}
  p_{11}(n) & p_{12}(n) & p_{13}(n) \\
  0 & p_{22}(n) & p_{23}(n) \\
  0 & 0 & p_{33}(n)
\end{pmatrix}
$$
where

$$
p_{11}(n)=\D\frac{(-1)^n(\alpha+\beta-k+n+2)(\alpha+\beta-k+n+3)(\alpha+1)_n}
{(\alpha+\beta-k+3)(\alpha+\beta-k+2)(\alpha+\beta+n+3)_{n}}
$$
$$
p_{12}(n)=\D\frac{2(-1)^nn(\beta-k+1)(\alpha+\beta-k+n+3)(\alpha+2)_{n-1}}
{(k+n+1)(\alpha+\beta-k+4)(\alpha+\beta-k+2)(\alpha+\beta+n+4)_{n-1}}
$$
$$
p_{13}(n)=\D\frac{(-1)^nn(n-1)(\beta-k+1)(\beta-k+2)(\alpha+3)_{n-2}}
{(k+n)(k+n+1)(\alpha+\beta-k+3)(\alpha+\beta-k+4)(\alpha+\beta+n+5)_{n-2}}
$$
$$
p_{22}(n)=\D\frac{(-1)^{n+1}(\alpha+\beta-k+2)(\alpha+\beta-k+n+4)(\alpha+2)_n}
{n(\alpha+\beta-k+4)(\alpha+\beta+n+4)_n}
$$
$$
p_{23}(n)=\D\frac{(-1)^{n+1}(\beta-k+2)(\alpha+\beta-k+2)(\alpha+3)_{n-1}}
{(k+n)(\alpha+\beta-k+4)(\alpha+\beta+n+5)_{n-1}}
$$
$$
p_{33}(n)=\D\frac{(-1)^n(\alpha+\beta-k+2)(\alpha+\beta-k+3)(\alpha+3)_n}
{n(n+1)(\alpha+\beta+n+5)_n}
$$
and $(a)_n$ denotes the shifted factorial defined by
$$
(a)_n=a(a+1)\cdots(a+n-1)\quad \mbox{for}\quad n>0,\;
(a)_0=1,\;(a)_{-1}=(a)_{-2}=0,
$$
with $a$ any real or complex number.

\subsection{The algebra of operators}\label{sec22}

For convenience we use in this section the conjugation of the second
order differential operator $D$ arising from representation theory
introduced in section \ref{sec1} instead of the one introduced in
subsection \ref{sec21}.

This will give rise to a new family of orthogonal polynomials which
are easily related to the ones obtained in the previous section.
This new family has certain computational advantages compared to the
previous one, but in principle one could use either one of them.

The resulting operator is:

\begin{align*}
D_1&=t(1-t)\D\frac{d^2}{dt^2}+ \left[\begin{pmatrix}
\alpha+3 & 0 & 0\\
-1 & \alpha+2 & 0\\
0 & -2 & \alpha+1
\end{pmatrix}-t \begin{pmatrix}
\alpha+\beta+4 & 0 & 0\\
0 & \alpha+\beta+5 & 0\\
0 & 0 & \alpha+\beta+6
\end{pmatrix}\right] \D\frac{d}{dt} \\ & + \begin{pmatrix}
0 & 2(\beta-k+1) & 0\\
0 & -(\alpha+\beta-k+2) & \beta-k+2\\
0 & 0 & -2(\alpha+\beta-k+3)
\end{pmatrix}.
\end{align*}

This operator has a sequence of (non--monic) matrix valued
orthogonal polynomials $\Pn$ with $P_0=I$ and
$$
D_1 P_n^*=P_n^*\Lambda_n(D_1),\;n=0,1,2,\ldots
$$
where the eigenvalue can be chosen (by an argument as the one given
in subsection \ref{sec21}) to be the same as in subsection
\ref{sec21}:
$$
\Lambda_n(D_1)= \begin{pmatrix}
t_1 & 0 & 0 \\
0 & t_2 & 0 \\
0 & 0 & t_3
\end{pmatrix}
$$
and the values of $t_i,\;i=1,2,3$ are given in (\ref{diag}).

An explicit expression of the coefficients of these $\Pn$, using
(\ref{coeff}), is given in the appendix in section \ref{app1}. These
expressions feature some denominators whose non--vanishing is
equivalent to the genericity assumptions made above.

The main goal of this section is to study the structure of the
following set:

$$
\mathcal{D}=\{D: DP_n^*=P_n^*\Lambda_n(D),\; n=0,1,2,\ldots\}.
$$

Here, $D$ is a differential operator of \emph{arbitrary order}.
Using the fact that the leading coefficient of $P_n^*(t)$ is a non
singular matrix we can see that $D=\sum_{j=0}^r F_j(t)\partial_t^j$,
$F_j(t)$ a matrix polynomial of degree $\leq j$, $j=0,\cdots, r$,
$r$ the order of $D$. $\mathcal{D}$ is clearly both a complex vector
space and an algebra under composition.

We observe that the map between differential operators and the
corresponding eigenvalues given by
$$
\Lambda_n:\mathcal{D}\longrightarrow M(3,\mathbb{C}),\quad
n=0,1,2,\ldots
$$
is a \emph{faithful representation}. The property
$\Lambda_n(D_1D_2)=\Lambda_n(D_1)\Lambda_n(D_2)$ with
$D_1,D_2\in\mathcal{D}$ is easy to show. If $\Lambda_n(D)=\Theta$
with $D\in\mathcal{D}$ for all $n$, then $DP_n^*=\Theta$ and we can
easily conclude that $D=\Theta$.

For all the calculations in the rest of this section it is very
useful to begin with monic matrix valued orthogonal polynomials and
then come back to the ones considered at the beginning. It is
important to note that the algebra $\mathcal{D}$ is independent of
the choice of the family of matrix valued orthogonal polynomials.
The eigenvalues are changed by an $n$ dependent conjugation as
noticed before.

We now set out to solve the equations $DP_n^*=P_n^*\Lambda_n(D)$
where $D=\sum_{j=0}^r F_j(t)\partial_t^j$, $F_j(t)$ a matrix
polynomial of degree $\leq j$, $j=0,\cdots, r$, $r$ the order of $D$
and $P_n^*$ as above.

In the following table, obtained by direct computations, we exhibit
the number of {\bf new} linearly independent differential operators
that appear as one increases the order of the operators in question:

\begin{center}\begin{tabular}{||c||c|c|c|c|c|c|c|c|c|c|c|c|c||}
  \hline order & 0 & 1 & 2 & 3 & 4 & 5 & 6 & 7 & 8 & 9 & 10 & 11 &
12\\\hline
  dimension & 1 & 0 & 2 & 0 & 3 & 0 & 3 & 0 & 3 & 0 & 3 & 0 & 3 \\ \hline
\end{tabular}\end{center}

There are no odd order differential operators, two linearly
independent second order differential operators and three {\bf new}
linearly independent differential operator in each even order
afterwards. We denote by $\mathfrak{D}_{2i}$ the space of
differential operators of order less or equal than $2i$ that have
our family of orthogonal polynomials as their eigenfunctions.

A possible second order differential operator linearly independent
from $D_1$ is the following:

\begin{eqnarray*}
D_2=\begin{pmatrix}
t(1-t) & 0 & 0\\
t/2 & t(1-t)/2 & 0\\
0 & -t & 0
\end{pmatrix}\D\frac{d^2}{dt^2}+\left[\begin{pmatrix}
\alpha+\beta-k+4 & \beta-k+1 & 0\\
-(\alpha+\beta-k+4)/2 & (\alpha+4)/2 & (\beta-k+2)/2\\
0 & -(\alpha+\beta-k+5) & -(\beta-k+2)
\end{pmatrix}\right.\bigskip\\ \left.-t \begin{pmatrix}
\alpha+\beta+4 & \beta-k+1 & 0\\
0 & (\alpha+\beta+5)/2 & (\beta-k+2)/2\\
0 & 0 & \alpha+\beta+6
\end{pmatrix}\right]\D\frac{d}{dt}+ \begin{pmatrix}
0 & -k(\beta-k+1) & 0\\
0 & k(\alpha+\beta-k+2)/2 & -k(\beta-k+2)/2\\
0 & 0 & k(\alpha+\beta-k+3)
\end{pmatrix}.
\end{eqnarray*}

The eigenvalue associated with $D_2$ is:
$$
\Lambda_n(D_2)= \begin{pmatrix}
-n^2-n(\alpha+\beta+3) & 0 & 0 \\
0 & \D\frac{1}{2}\bigl(- n^{2}
-n(\alpha+\beta+4)+k(\alpha+\beta-k+2)\bigr) & 0
\\
0 & 0 & k(\alpha+\beta-k+3)
\end{pmatrix}.
$$

This means in particular that $[D_1,D_2]=D_1D_2-D_2D_1=\Theta$.

We come now to the problem of exhibiting, for each value of
$i,\;i=2,3,\ldots$ a set of three linearly independent differential
operators that will have order $2i$ and will span a subspace
$\mathfrak{R}_{2i}$ such that
$$\mathfrak{D}_{2i}=\mathfrak{D}_{2i-2}\oplus \mathfrak{R}_{2i}.$$

In the case of $i=2$ we can choose the set $\mathfrak{R}_4$
consisting of $\{D_1^2,D_2^2,D_1D_2\}$. One can easily check that
these operators are linearly independent. Moreover, it is also easy
to check, using the isomorphic bijection between differential
operators and the corresponding eigenvalues, that the set
$\{D_1^i,D_2^i,D_1^{i-1}D_2\}$, $i=2,3,\ldots$, is linearly
independent and can be chosen as a basis for $\mathfrak{R}_{2i}$.
Since all these operators are obtained from $D_1$ and $D_2$, we can
\textbf{conjecture} that the algebra $\mathcal{D}$ is \emph{abelian}
and coincides with the subalgebra $\mathcal{A}$ generated by
$\{I,D_1,D_2\}$.

We turn our attention now to finding relations among these
generators.

One relation follows from what we said earlier regarding
$\mathfrak{R}_6$. The operator $D_1D_2^2$ cannot be linearly
independent from the elements in the basis of $\mathfrak{R}_6$. One
is thus inclined to look for $s_1, s_2,\ldots,s_9$ such that:

\begin{equation}\label{coef1}
s_1I+s_2D_1+s_3D_2+s_4D_1^2+s_5D_2^2+s_6D_1D_2+s_7D_1^3+s_8D_2^3+s_9D_1^2D_2=D_1D_2^2.
\end{equation}

The result is apparently rather unilluminating, namely we get

\begin{equation}\label{coef2}
   \begin{split}
      s_1&=0\\
      s_2&=-\frac{1}{3}k(k+1)(\alpha+\beta-k+3)(\alpha+\beta-k+2)\\
      s_3&=\frac{1}{3}k(k+1)(\alpha+\beta-k+3)(\alpha+\beta-k+2)\\
      s_4&=-\frac{1}{3}k(\alpha+\beta-k+3)\\
s_5&=-\frac{2}{3}-\frac{1}{3}\alpha-k(\beta+1)+k^2-k\alpha-\frac{1}{3}\beta-\frac{4}{3}k\\
s_6&=\frac{2}{3}+\frac{1}{3}\alpha+\frac{4}{3}k(\beta+1)-\frac{4}{3}k^2+\frac{4}{3}k\alpha+\frac{1}{3}\beta+2k\\
      s_7&=0,\quad s_8=\frac{2}{3},\quad s_9=\frac{1}{3}.
      \end{split}
    \end{equation}

Now something remarkable happens: if we put (\ref{coef2}) in
(\ref{coef1}), one can obtain the following amazing
\emph{factorization}:

\begin{equation}\label{fact}
\bigl(D_1-D_2\bigr)\bigl(D_2-k(\alpha+\beta-k+3)\bigr)\bigl(D_1-2D_2+(1+k)(\alpha+\beta-k+2)\bigr)=\Theta.
\end{equation}

We are thankful to an anonymous referee for pointing out a way to
make this factorization quite transparent. The argument goes as
follows:

Since the algebra of operators generated by $D_1$ and $D_2$ is
isomorphic to the quotient of $\mathbb{C}[x,y]$ (the algebra of
polynomials in two commuting variables) by the ideal of polynomials
$p\in\mathbb{C}[x,y]$ that satisfy $p\,(D_1,D_2)=\Theta$ and the
mapping $ \Lambda_n:\mathcal{D}\longrightarrow M(3,\mathbb{C})$ is a
faithful representation as indicated earlier, it follows that the
condition $p\,(D_1,D_2)=\Theta$ is equivalent to
$$
\Lambda_n(p(D_1,D_2))= p
(\Lambda_n(D_1),\Lambda_n(D_2))=\Theta,\quad n=0,1,2,\cdots.
$$

Since the matrices $\Lambda_n(D_1)$ and $\Lambda_n(D_2)$  are
diagonal with eigenvalues $t_1(n)$, $t_2(n)$ , $t_3(n)$ and $r_1(n)=
t_1(n)$, $r_2(n)= \frac{1}{2}(t_2(n)+(k+1)(\alpha+\beta-k+2))$,
$r_3(n)=k(\alpha+\beta-k+3)$, respectively,  it is easy to spot
three polynomials that should be included as factors of any
polynomial in the ideal in question. These polynomials are
\begin{equation*}
   \begin{split}
      x-y&\quad \mbox{from the entry}\quad (1,1)\\
      x-2y +(1+k)(\alpha+\beta-k+2)&\quad \mbox{from the entry} \quad(2,2) \quad \mbox{and}\\
      y-k(\alpha+\beta-k+3)&\quad \mbox{from the entry} \quad(3,3).
      \end{split}
    \end{equation*}
These are the factors that appear in the simplest of the relations
that is obtained by multiplying these factors and replacing $x,y$ by
$D_1,D_2$, namely (\ref{fact}).

One can consider this {\bf curve} as the analog of the
Burchnall-Chaundy curve, an algebraic geometric object associated to
a commutative algebra of differential operators acting on scalar
functions. For a sample of the literature on this beautiful and deep
subject which traces its origin to some early work of J.Burchnall
and T. Chaundy, see \cite {I} and its references. This has become
more recently a very active area and one can consult
\cite{Kr,Mum,vMM}.

It is worth noting that in the scalar case, as long as the
eigenfunctions are orthogonal polynomials and the differential
operators are of order two, the algebra is trivial. See \cite{M1}
for a proof.

It is important to note that it is only the product of the three
factors in (\ref{fact}), and not the product of any smaller subset
of them, that vanishes.

\section{The two steps example}

\subsection{Generating the polynomial eigenfunctions of the differential
operator}\label{sec31}

We consider now the example resulting from the {\em two steps
situation} described in the second half of section \ref{sec1}. This
subsection proceeds along the lines of subsection \ref{sec21}.

As in the one step case, we retain the symbol $D$ for the new
differential operator given in {\em hypergeometric form} by the
expression:

$$
D=t(1-t)\D\frac{d^2}{dt^2}+(X-tU)\D\frac{d}{dt}+V
$$
where
$$
V= \begin{pmatrix}
0 & 0 & 0 & 0\\
0 & -(\alpha+\beta-k_1+3) & 0 & 0\\
0 & 0 & -(\alpha+\beta-k_2+2) & 0 \\
0 & 0 & 0 & -(2\alpha+2\beta-k_1-k_2+6)
\end{pmatrix}
$$
$$
U= \begin{pmatrix} \alpha+\beta+4 & -1 & -1
&\frac{2\alpha+2\beta-k_1-k_2+6}{(\alpha+\beta-k_2+2)(\alpha+\beta-k_1+3)}\\
0 & \alpha+\beta+5 & 0 &
-\frac{(k_1-k_2)(\alpha+\beta-k_1+4)}{(\alpha+\beta-k_1+3)(k_1-k_2-1)}\\
0 & 0 & \alpha+\beta+5 &
-\frac{(k_1-k_2-2)(\alpha+\beta-k_2+3)}{(k_1-k_2-1)(\alpha+\beta-k_2+2)}\\
0 & 0 & 0 &\alpha+\beta+6
\end{pmatrix}
$$

$$
X= \begin{pmatrix} C_{11}(\alpha,\beta,k_1,k_2)&
\frac{(\alpha+1)(\alpha+\beta-k_2+3)}{(\alpha+\beta-k_2+2)(\alpha+\beta-k_1+3)}
&
\frac{(\alpha+1)(\alpha+\beta-k_1+4)}{(\alpha+\beta-k_2+2)(\alpha+\beta-k_1+3)}&
0\\
\frac{(k_1-k_2)(\beta-k_1+2)}{(\alpha+\beta-k_1+3)(k_1-k_2-1)} &
C_{22}(\alpha,\beta,k_1,k_2) &C_{23}(\alpha,\beta,k_1,k_2)
& \frac{(k_1-k_2)(\alpha+2)}{(\alpha+\beta-k_2+3)(k_1-k_2-1)}\\
\frac{(k_1-k_2-2)(\beta-k_2+1)}{(k_1-k_2-1)(\alpha+\beta-k_2+2)} &
C_{32}(\alpha,\beta,k_1,k_2) & C_{33}(\alpha,\beta,k_1,k_2) &
\frac{(k_1-k_2-2)(\alpha+2)}{(\alpha+\beta-k_1+4)(k_1-k_2-1)}\\
0 &
\frac{(\alpha+\beta-k_1+3)(\beta-k_2+1)}{(\alpha+\beta-k_1+4)(\alpha+\beta-k_2+3)}
&
\frac{(\alpha+\beta-k_2+2)(\beta-k_1+2)}{(\alpha+\beta-k_1+4)(\alpha+\beta-k_2+3)}
& C_{44}(\alpha,\beta,k_1,k_2)
\end{pmatrix}
$$
where
$$
C_{11}(\alpha,\beta,k_1,k_2)=\frac{(\alpha+1)(\alpha+\beta-k_1+4)(\alpha+\beta-k_2+3)}
{(\alpha+\beta-k_2+2)(\alpha+\beta-k_1+3)}
$$
$$
C_{23}(\alpha,\beta,k_1,k_2)=\frac{(k_1-k_2)(\beta-k_1+2)}{(k_1-k_2-1)(\alpha+\beta-k_1+3)(\alpha+\beta-k_2+3)}
$$
$$
C_{32}(\alpha,\beta,k_1,k_2)=\frac{(k_1-k_2-2)(\beta-k_2+1)}{(k_1-k_2-1)(\alpha+\beta-k_1+4)(\alpha+\beta-k_2+2)}.
$$
The elements $C_{22}, C_{33}, C_{44}$ are given in the appendix in
section \ref{app2}.

A sequence of (non--monic) matrix valued orthogonal polynomials
$\Pn$ with $P_0=I$ is then obtained by solving the following
differential equation:

$$
t(1-t)\D\frac{d^2}{dt^2}P_n^*+(X-tU)\D\frac{d}{dt}P_n^*+VP_n^*=P_n^*\Lambda_n,
$$
where

$$
\Lambda_n= \begin{pmatrix}
t_1 & 0 & 0 & 0 \\
0 & t_2 & 0 & 0 \\
0 & 0 & t_3 & 0\\
0 & 0 & 0 & t_4
\end{pmatrix}
$$
and $t_i$, $i=1,\ldots,4$ are given below in (\ref{diag4}). The
reasoning behind this choice is exactly the same as in subsection
\ref{sec21} and it is not repeated here.
\par\bigskip

Proceeding as in the previous section along the lines in \cite{T2},
we replace the differential operator by another operator acting on
functions that take values on $\mathbb{C}^{16}$ and consider the
left and right multiplication by matrices in $M(4,\mathbb{C})$ as
linear maps in $\mathbb{C}^{16}$. Thus we can consider the following
equivalent differential equation:

$$
t(1-t)\D\frac{d^2}{dt^2}\mbox{vec}(P_n^*)+(C-t\tilde{U})\D\frac{d}{dt}\mbox{vec}(P_n^*)-\tilde{T}\mbox{vec}(P_n^*)=\Theta
$$
where $C$ and $\widetilde{U}$ are the $16\times16$ matrices obtained
from $X$ and $U$ respectively, in the same manner as was done
earlier. A similar procedure is applied to obtain the matrix
$\widetilde{T}$ and we get:

$$
\tilde{T}=\mbox{diag}\{t_1,t_2,t_3,t_4,t_5,t_6,t_7,t_8,t_9,t_{10},t_{11},t_{12},t_{13},t_{14},t_{15},t_{16}\}
$$
where
\begin{equation}\label{diag4}
   \begin{split}
      t_1&=-n^2-n(\alpha+\beta+3)\\
      t_2&=- n^{2} -n(\alpha+\beta+4)-(\alpha+\beta-k_1+3)\\
      t_3&=-n^{2} -n(\alpha+\beta+4)-(\alpha+\beta-k_2+2)\\
      t_4&=- n^{2}-n(\alpha+\beta+5)-(2\alpha+2\beta-k_1-k_2+6)\\
      t_5&=-n^2-n(\alpha+\beta+3)+\alpha+\beta-k_1+3\\
      t_6&=- n^{2} -n(\alpha+\beta+4)\\
      t_7&=- n^{2} -n(\alpha+\beta+4)+k_2-k_1+1\\
      t_8&=- n^{2} -n(\alpha+\beta+5)-(\alpha+\beta-k_2+3)\\
      t_9&=-n^2-n(\alpha+\beta+3)+\alpha+\beta-k_2+2\\
      t_{10}&=- n^{2} -n(\alpha+\beta+4)+k_1-k_2-1\\
      t_{11}&=- n^{2} -n(\alpha+\beta+4)\\
      t_{12}&=- n^{2}-n(\alpha+\beta+5)-(\alpha+\beta-k_1+4)\\
      t_{13}&=-n^2-n(\alpha+\beta+3)+2\alpha+2\beta-k_1-k_2+6\\
      t_{14}&=- n^{2} -n(\alpha+\beta+4)+\alpha+\beta-k_2+3\\
      t_{15}&=- n^{2} -n(\alpha+\beta+4)+\alpha+\beta-k_1+4\\
      t_{16}&=- n^{2} -n(\alpha+\beta+5)
    \end{split}
    \end{equation}

Now we find a pair of matrices $A,B$ as in the previous section,
satisfying the non--linear matrix equations
$$\tilde{U}=I+A+B\quad \mbox{and}\quad
\tilde{T}=AB.$$ As in the previous section, we pick

$$
A=\begin{pmatrix}
  A_{11} & A_{12} & A_{13} & A_{14} \\
  \Theta & A_{22} & \Theta & A_{24} \\
  \Theta & \Theta & A_{33} & A_{34} \\
  \Theta & \Theta & \Theta & A_{44}
\end{pmatrix}
$$
$$
B=\tilde{U}-A-I
$$
where each $A_{ij}$ is a $4\times 4$ diagonal matrix.

Observe that $A_{23}=\Theta$, a fact which will have important
consequences.

We put

$$
A_{ii}=\begin{pmatrix}
  \phi_i & 0 & 0 & 0 \\
  0 & \varphi_i & 0 & 0 \\
  0 & 0 & \psi_i & 0 \\
  0 & 0 & 0 & \xi_i
\end{pmatrix}
\;\; i=1,2,3,4.
$$
with $\phi_i,\varphi_i,\psi_i,\xi_i\;i=1,2,3,4$ to be determined
later.

The elements of the diagonal of $A_{12}$ are of the form:

$$
\D\frac{\gamma_1}{-\gamma_1-\gamma_2+\alpha+\beta+4},
$$
with $\gamma=\phi,\varphi,\psi,\xi$ in each entry.

The elements of the diagonal of $A_{13}$ are of the form:

$$
\D\frac{\gamma_1}{-\gamma_1-\gamma_3+\alpha+\beta+4},
$$
with $\gamma=\phi,\varphi,\psi,\xi$ in each entry.

If we denote by

$$
\omega_{14}=\D\frac{2\alpha+2\beta-k_1-k_2+6}{(\alpha+\beta-k_2+2)(\alpha+\beta-k_1+3)}
$$
$$
\omega_{24}=-\D\frac{(k_1-k_2)(\alpha+\beta-k_1+4)}{(\alpha+\beta-k_1+3)(k_1-k_2-1)}
$$
$$
\omega_{34}=-\D\frac{(k_1-k_2-2)(\alpha+\beta-k_2+3)}{(k_1-k_2-1)(\alpha+\beta-k_2+2)}
$$
the elements of the last column of $U$, the elements of the diagonal
of $A_{34}$ are of the form:

$$
\D\frac{-\omega_{34}\gamma_3}{-\gamma_3-\gamma_4+\alpha+\beta+5},
$$
with $\gamma=\phi,\varphi,\psi,\xi$ in each entry, and the elements
of the diagonal of $A_{24}$ are of the form:

$$
\D\frac{-\omega_{24}\gamma_2}{-\gamma_2-\gamma_4+\alpha+\beta+5},
$$
with $\gamma=\phi,\varphi,\psi,\xi$ in each entry.

The elements of the diagonal of $A_{14}$ are a bit more complicated,
they are given by:

\begin{eqnarray*}
-\frac{\omega_{14}\gamma_1}{-\gamma_1-\gamma_4+\alpha+\beta+5}-\frac{\omega_{24}\gamma_1}{(-\gamma_1-\gamma_4+\alpha+\beta+5)
(-\gamma_1-\gamma_2+\alpha+\beta+4)}
\bigskip\\-\frac{\omega_{34}\gamma_1}{(-\gamma_1-\gamma_4+\alpha+\beta+5)(-\gamma_1-\gamma_3+\alpha+\beta+4)}\bigskip\\
-\frac{\omega_{24}\gamma_1\gamma_2}{(-\gamma_1-\gamma_4+\alpha+\beta+5)(-\gamma_1-\gamma_2+\alpha+\beta+4)
(-\gamma_2-\gamma_4+\alpha+\beta+5)}\bigskip\\
-\frac{\omega_{34}\gamma_1\gamma_3}{(-\gamma_1-\gamma_4+\alpha+\beta+5)(-\gamma_1-\gamma_3+\alpha+\beta+4)(-\gamma_3-\gamma_4+\alpha+\beta+5)},
\end{eqnarray*}

with $\gamma=\phi,\varphi,\psi,\xi$ in each entry.

The parameters $\phi_i,\varphi_i,\psi_i,\xi_i$, $i=1,2,3,4$ are
subject to the following conditions:\\
$$
\left\{\begin{array}{l}
\phi_1=-n\quad \mbox{or}\quad \phi_1=n+\alpha+\beta+3\\
\varphi_1^2-(\alpha+\beta+3)\varphi_1+t_2=0\\
\psi_1^2-(\alpha+\beta+3)\psi_1+t_3=0 \\
\xi_1^2-(\alpha+\beta+3)\xi_1+t_4=0
\end{array}\right.
$$
$$
\left\{\begin{array}{l}
\phi_2^2-(\alpha+\beta+4)\phi_2+t_5=0\\
\varphi_2=-n\quad \mbox{or}\quad \varphi_2=n+\alpha+\beta+4\\
\psi_2^2-(\alpha+\beta+4)\psi_2+t_7=0\\
\xi_2^2-(\alpha+\beta+4)\xi_2+t_8=0
\end{array}\right.
$$
$$
\left\{\begin{array}{l}
\phi_3^2-(\alpha+\beta+4)\phi_2+t_9=0\\
\varphi_3^2-(\alpha+\beta+4)\varphi_3+t_{10}=0\\
\psi_3=-n\quad \mbox{or}\quad \psi_3=n+\alpha+\beta+4\\
\xi_3^2-(\alpha+\beta+4)\xi_3+t_{12}=0
\end{array}\right.
$$
$$
\left\{\begin{array}{l}
\phi_4^2-(\alpha+\beta+5)\phi_4+t_{13}=0\\
\varphi_4^2-(\alpha+\beta+5)\varphi_4+t_{14}=0\\
\psi_4^2-(\alpha+\beta+5)\psi_4+t_{15}=0\\
\xi_4=-n\quad \mbox{or}\quad \xi_4=n+\alpha+\beta+5
\end{array}\right.
$$

Now, $\mbox{vec}(P_n^*)$ satisfies the matrix hypergeometric
equation:

\begin{equation*}
t(1-t)\D\frac{d^2}{dt^2}\mbox{vec}(P_n^*)+(C-t(I+A+B))\D\frac{d}{dt}\mbox{vec}(P_n^*)-AB\mbox{vec}(P_n^*)=\Theta.
\end{equation*}

The eigenvalues of $C$ are $\{\alpha+1,\alpha+2,\alpha+3\}$ with
$\alpha+1$ and $\alpha+3$ with multiplicity 4 and $\alpha+2$ with
multiplicity 8. So if $\alpha>-1$ we meet the conditions, see
\cite{T2}, of the definition of the \emph{hypergeometric function},
as in the one step case.

The function $_2F_1(C,A,B;t)$ is not a polynomial function, as in
the classical case, but nevertheless we get, as in \cite{T2}, the
family of vector valued polynomials of degree $n$ in $t$:

$$
\mbox{vec}\bigl(P_n^*(t)\bigr)=\,_2F_1(C,A,B;t)\mbox{vec}\bigl(P_n^*(0)\bigr).
$$

In this case we can also give explicitly the value of $P_n^*(0)$, as
before, using (\ref{coeff}):

$$
P_n^*(0)=\begin{pmatrix}
  p_{11}(n) & p_{12}(n) & p_{13}(n) & p_{14}(n) \\
  0 & p_{22}(n) & 0 & p_{24}(n) \\
  0 & 0 & p_{33}(n) & p_{34}(n) \\
  0 & 0 & 0 & p_{44}(n)
\end{pmatrix},
$$
where

$$
p_{11}(n)=\D\frac{(-1)^n(\alpha+\beta-k_1+n+3)(\alpha+\beta-k_2+n+2)(\alpha+1)_n}
{(\alpha+\beta-k_1+3)(\alpha+\beta-k_2+2)(\alpha+\beta+n+3)_n}
$$
$$
p_{12}(n)=\D\frac{(-1)^nn(k_2-k_1)(\beta-k_1+2)(\alpha+\beta-k_2+n+2)(\alpha+2)_{n-1}}
{(k_1+n)(k_2-k_1+1)(\alpha+\beta-k_1+3)(\alpha+\beta-k_2+3)(\alpha+\beta+n+4)_{n-1}}
$$
$$
p_{13}(n)=\D\frac{(-1)^nn(k_2-k_1+2)(\beta-k_2+1)(\alpha+\beta-k_1+n+3)(\alpha+2)_{n-1}}
{(k_2+n+1)(k_2-k_1+1)(\alpha+\beta-k_1+4)(\alpha+\beta-k_2+2)(\alpha+\beta+n+4)_{n-1}}
$$
$$p_{14}(n)=\D\frac{(-1)^nn(n-1(\beta-k_2+1)(\beta-k_1+2)(\alpha+3)_{n-2}}
{(k_1+n)(k_2+n+1)(\alpha+\beta-k_1+4)(\alpha+\beta-k_2+3)(\alpha+\beta+n+5)_{n-2}}
$$
$$
p_{22}(n)=\D\frac{(-1)^{n+1}(\alpha+\beta-k_2+n+4)(\alpha+\beta-k_1+3)(\alpha+2)_n}
{n(\alpha+\beta-k_2+3)(\alpha+\beta+n+4)_n}
$$
$$
p_{24}(n)=\D\frac{(-1)^{n+1}(\beta-k_2+1)(\alpha+\beta-k_1+3)(\alpha+3)_{n-1}}
{(k_2+n+1)(\alpha+\beta-k_2+3)(\alpha+\beta+n+5)_{n-1}}
$$
$$
p_{33}(n)=\D\frac{(-1)^{n+1}(\alpha+\beta-k_1+n+4)(\alpha+\beta-k_2+2)(\alpha+2)_n}
{(\alpha+\beta-k_1+4)(\alpha+\beta+n+4)_n}
$$
$$
p_{34}(n)=\D\frac{(-1)^{n+1}(\beta-k_1+2)(\alpha+\beta-k_2+2)(\alpha+3)_{n-1}}
{(k_1+n)(\alpha+\beta-k_1+4)(\alpha+\beta+n+5)_{n-1}}
$$
$$
p_{44}(n)=\D\frac{(-1)^n(\alpha+\beta-k_1+3)(\alpha+\beta-k_2+2)(\alpha+3)_n}
{n(n+1)(\alpha+\beta+n+5)_n}.
$$

Note that $p_{23}(n)=0$. This is related to the vanishing of the
$(2,3)$--block of $A$.

\subsection{The algebra of operators}\label{sec32}

In this case, just as in the one step example, we use a different
conjugation of the second order differential operator arising from
representation theory and introduced at the end of section
\ref{sec1}. This is done only for computational convenience, and all
these operators (and the corresponding orthogonal polynomials) are
related to each other by simple conjugations.

The operator now is:
\begin{equation*}
   D_1= A_2(t)\frac{d^2}{dt^2}+
A_1(t)\frac{d}{dt}+ A_0(t),
\end{equation*}
with $ A_2$, $A_1$, $A_0$ given by:
\begin{eqnarray*}
A_2(t) &= &t(1-t) \\
A_1(t) &=  &\begin{pmatrix}
\alpha+3 &0 &0 &0\\
-1 & \alpha+2 &0 &0 \\
-1 &0 & \alpha+2 &0\\
0 & -\tfrac{k_2-k_1+2}{k_2-k_1+1} & -\tfrac{k_2-k_1}{k_2-k_1+1}
&\alpha+1
\end{pmatrix} \\
& &\quad - t \begin{pmatrix}
\alpha+\beta+4 & 0&0 &0 \\
0 &\alpha+\beta+5 &0 &0\\
0& 0& \alpha+\beta+5\\
0&0 &0 &\alpha+\beta+6
\end{pmatrix}  \\
A_0(t) &= &\begin{pmatrix} 0 &
\tfrac{(k_2-k_1+2)(\beta-k_2+1)}{k_2-k_1+1} &
\tfrac{(k_2-k_1)(\beta-k_1+2)}{k_2-k_1+1}&0 \\
0 &-(\alpha+\beta+2)+k_2 &0 &\beta-k_1+2 \\
0&0& -(\alpha+\beta+3)+k_1 & \beta-k_2+1\\
0&0 &0 &-2(\alpha+\beta+3)+k_1+k_2
\end{pmatrix}.
\end{eqnarray*}
The eigenvalue associated with $D_1$ can be chosen, as before, in
the form:

$$
\Lambda_n(D_1)=\left(%
\begin{array}{cccc}
  t_1 & 0 & 0 & 0 \\
  0 & t_3& 0 & 0 \\
  0 & 0& t_2 & 0 \\
  0 & 0 & 0 & t_4 \\
\end{array}%
\right)
$$
where $t_i, i=1,2,3,4$ are given in (\ref{diag4}). Note that $t_2$
and $t_3$ above do not appear now in the natural order.

An explicit expression of the coefficients of these $\Pn$, using
(\ref{coeff}), is given in the appendix in section \ref{app3}. These
expressions feature some denominators whose non--vanishing is
equivalent to the genericity assumptions implicitly made in
subsection \ref{sec31}.

We now consider the analog of the algebra $\mathcal{D}$ from
subsection \ref{sec22}. Recall that these differential operators
have matrix valued polynomial coefficients of the appropriate
degrees.

One finds that the commutator of the operator $D_1$, within the
class of operators of order less or equal than two, has dimension
\emph{four}. A basis is made up of three operators of order two and
the identity. One of the linearly independent operators is given by:

\begin{equation*}
   D_2= B_2(t)\frac{d^2}{dt^2}+ B_1(t)\frac{d}{dt}+ B_0(t),
\end{equation*}
with $ B_2$, $B_1$, $B_0$ given by:
\begin{eqnarray*}
B_2(t) & = & \begin{pmatrix} 0 & 0 & 0 & 0 \\
0&0&0&0\\
\frac{k_1-k_2-1}{k_1-k_2}t &0 & \frac{k_1-k_2-1}{k_1-k_2}t(1-t) & 0 \\
0& \frac{k_1-k_2-2}{k_1-k_2-1}t & \frac{1}{k_1-k_2-1}t & t(1-t)
\end{pmatrix}\\
B_1(t) & =  & \begin{pmatrix}  -(\beta-k_1+2)&
 0&  -(\beta-k_1+2) & 0 \\
0  & -\frac{(\beta-k_1+2)(k_1-k_2-2)}{k_1-k_2-1} &
-\frac{\beta-k_1+2}{k_1-k_2-1} &
 -(\beta-k_1+2) \\
\frac{(\alpha+\beta-k_1+3)(k_1-k_2-1)}{k_1-k_2} &
-\frac{(\beta-k_2+1)(k_1-k_2-2)}{(k_1-k_2-1)(k_1-k_2)} &
\frac{b_{33}}{(k_1-k_2-1)(k_1-k_2)} &
-\frac{\beta-k_2+1}{k_1-k_2}\\
0 & \frac{(k_1-k_2-2)(\alpha+\beta-k_1+4)}{k_1-k_2-1}  &
\frac{\alpha+\beta-k_1+4}{k_1-k_2-1} &
\alpha+\beta-k_1+4 \\
\end{pmatrix} \\
& &\quad - t \begin{pmatrix} 0 &
 0&
 -(\beta-k_1+2)& 0 \\
0& 0 & 0&
  -(\beta-k_1+2) \\
0 & 0 & \frac{(\alpha+\beta+5)(k_1-k_2-1)}{k_1-k_2}&
-\frac{\beta-k_2+1}{k_1-k_2} \\ 0& 0& 0& \alpha+\beta+6 \\
\end{pmatrix}  \\
B_0(t) & = &\begin{pmatrix} 0&0 &
(1+k_1)(\beta-k_1+2)& 0 \\
0 &0 &
0&  (1+k_1)(\beta-k_1+2)\\
0& 0& -\frac{(1+k_1)(\alpha+\beta-k_1+3)(k_1-k_2-1)}{k_1-k_2}&
\frac{(1+k_1)(\beta-k_2+1)}{k_1-k_2} \\
0&0&0&-(1+k_1)(\alpha+\beta-k_1+4)
\end{pmatrix},
\end{eqnarray*}
where
\begin{eqnarray*}
\lefteqn{b_{33}=2-2\alpha k_2k_1-2\alpha
k_1-2k_1k_2\beta-8k_1k_2-k_1^3+\alpha+2k_1^2k_2+2k_2\beta+k_1^2\beta}\\
\noalign{\smallskip} & &\qquad\qquad +k_2^2\beta-2k_1\beta+2\alpha
k_2+\alpha k_2^2+\alpha k_1^2+7k_2-7k_1+3k_2^2+5k_1^2-k_1k_2^2.
\end{eqnarray*}
The eigenvalue associated with $D_2$ is:

\begin{equation*}
\Lambda_n(D_2)=\left(%
\begin{array}{cccc}
  0 & 0 & 0 & 0 \\
  0 & 0& 0 & 0 \\
  0 & 0& -\frac{(n+k_1+1)(n+\alpha+\beta-k_1+3)(k_1-k_2-1)}{k_1-k_2} & 0 \\
  0 & 0 & 0 &   -(n+k_1+1)(n+\alpha+\beta-k_1+4)
\end{array}%
\right).
\end{equation*}

The last second order linearly independent operator is given by:

\begin{equation*}
   D_3= C_2(t)\frac{d^2}{dt^2}+ C_1(t)\frac{d}{dt}+ C_0(t),
\end{equation*}
with $ C_2$, $C_1$, $C_0$ given by:
\begin{eqnarray*}
C_2(t) & = & \begin{pmatrix} 0 & 0 & 0 & 0 \\
\frac{1}{k_1-k_2-2}t&\frac{1}{k_1-k_2-2}t(1-t)&0&0\\
-t &0 & -t(1-t) & 0 \\
0& -t & 0 & -t(1-t)
\end{pmatrix}\\
C_1(t) & =  & \begin{pmatrix}  \beta-k_1+1&
 \frac{\beta-k_2+1}{k_1-k_2-1}&  \frac{(k_1-k_2)(\beta-k_1+2)}{k_1-k_2+1} & 0 \\
\frac{\alpha+\beta-k_2+2}{k_1-k_2-2}  &
\frac{c_{22}}{(k_1-k_2-2)(k_1-k_2-1)} &
\frac{(\beta-k_1+2)(k_1-k_2)}{(k_1-k_2-2)(k_1-k_2-1)} &
 \frac{(\beta-k_1+2)(k_1-k_2-1)}{k_1-k_2-2} \\
-(\alpha+\beta-k_1+3) & \frac{\beta-k_2+1}{k_1-k_2-1} &
-\frac{c_{33}}{k_1-k_2-1} &
0\\
0 & -\frac{c_{33}-\beta+k_1-1}{k_1-k_2-1}  &
\frac{k_1-k_2}{k_1-k_2-1} &
-(\alpha+\beta-k_1+3) \\
\end{pmatrix} \\
& &\quad - t \begin{pmatrix} 0 &
 -\frac{\beta-k_2+1}{k_1-k_2-1}&
 \frac{(k_1-k_2)(\beta-k_1+2)}{k_1-k_2-1}& 0 \\
0& \frac{\alpha+\beta+5}{k_1-k_2-2} & 0&
  \frac{(\beta-k_1+2)(k_1-k_2-1)}{k_1-k_2-2} \\
0 & 0 & -(\alpha+\beta+5)&
0 \\ 0& 0& 0& -(\alpha+\beta+6) \\
\end{pmatrix}  \\
C_0(t) & = &\begin{pmatrix} 0&\frac{(2+k_2)(\beta-k_2+1)}{k_1-k_2-1}
&
-\frac{(1+k_1)(k_1-k_2)(\beta-k_1+2)}{k_1-k_2-1}& 0 \\
0 &-\frac{(\alpha+\beta-k_2+2)(2+k_2)}{k_1-k_2-2} &
0&  \frac{(\beta-k_1+2)(k_2+2-k_1^2+k_1k_2)}{k_1-k_2-2}\\
0& 0& \mbox{{\footnotesize $(1+k_1)(\alpha+\beta-k_1+3)$}}&
-(\beta-k_2+1) \\
0&0&0&c_{44}
\end{pmatrix},
\end{eqnarray*}
where
\begin{eqnarray*}
\lefteqn{c_{22}=\alpha
k_1-2k_1\beta-5k_1+2k_2\beta+5k_2-8k_1k_2-\alpha
k_2+k_1^2\beta}\\
\noalign{\smallskip} & &\quad
+5k_1^2+k_2^2\beta+3k_2^2-k_1^3+2k_1^2k_2-2k_1k_2\beta-\alpha-k_1k_2^2
\end{eqnarray*}
\begin{equation*}
c_{33}=-\alpha k_2-4k_2-k_2\beta+k_1k_2-2+\alpha
k_1+k_1\beta+4k_1-\alpha-k_1^2
\end{equation*}
\begin{equation*}
c_{44}=-k_2+2\alpha+2k_1+\alpha k_1+\beta k_1-k_1^2+4+2n.
\end{equation*}
The eigenvalue associated with $D_3$ is:

\begin{equation*}
\Lambda_n(D_3)=\left(%
\begin{array}{cccc}
  0 & 0 & 0 & 0 \\
  0 & -\frac{(n+k_2+2)(n+\alpha+\beta-k_2+2)}{k_1-k_2-2}& 0 & 0 \\
  0 & 0& (n+k_1+1)(n+\alpha+\beta-k_1+3) & 0 \\
  0 & 0 & 0 &  \lambda_n\\
\end{array}%
\right),
\end{equation*}
where
$$\lambda_n=n^2+n(\alpha+\beta+5)+k_1(\alpha+\beta-k_1+2)
+2\alpha+2\beta-k_2+6.$$

We show in the following table, by direct calculations, the growth
of the dimensions of the space of differential operators in
$\mathcal{D}$ as one increases the order:

\begin{center}\begin{tabular}{||c||c|c|c|c|c|c|c|c|c|c|c|c|c||}
  \hline order & 0 & 1 & 2 & 3 & 4 & 5 & 6 & 7 & 8 & 9 & 10 & 11 &12\\\hline
  dimension & 1 & 0 & 3 & 0 & 6 & 0 & 6 & 0 & 6 & 0 & 6 & 0& 6\\ \hline
\end{tabular}\end{center}

There are no odd order differential operators, three linearly
independent second order differential operators, given by $D_1$,
$D_2$ and $D_3$ and six {\bf new} linearly independent differential
operator in each even order afterwards. We denote by
$\mathfrak{D}_{2i}$ the space of differential operators of order
less or equal than $2i$ that have our family of orthogonal
polynomials as their eigenfunctions.

We come now to the problem of exhibiting, for each value of
$i,\;i=2,3,\ldots$ a set of six linearly independent differential
operators that will have order $2i$ and will span a subspace
$\mathfrak{R}_{2i}$ such that
$$\mathfrak{D}_{2i}=\mathfrak{D}_{2i-2}\oplus \mathfrak{R}_{2i}.$$

For each $i,\;i=2,3,\ldots$, $\mathfrak{R}_{2i}$ has a basis made up
of \textbf{four} differential operators that commute with each other
and \textbf{two} extra ones that do not commute with any other
element in the base. This decomposition of $\mathfrak{R}_{2i}$ is
spelled out in the table below.

We denote the two extra operators in $\mathfrak{R}_{4}$ alluded to
above by $E$ and $F$.

The corresponding eigenvalues are given by:

$$
\Lambda_n(E)=(n+\alpha+\beta-k_1+3)(n+\alpha+\beta-k_2+3)(n+k_1)(n+k_1+1)E_{23}
$$
$$
\Lambda_n(F)=(n+\alpha+\beta-k_1+4)(n+\alpha+\beta-k_2+2)(n+k_2+1)(n+k_2+2)E_{32},
$$
where $E_{ij}$ denotes the $4\times4$ matrix with entry $(i,j)$
equal 1 and 0 otherwise.

In the following table we summarize the decomposition of
$\mathfrak{R}_{2i}$ discussed above:
\begin{center}
\begin{tabular}{c|c|c|}

   & COMM & NON--COMM \\\hline
  $\mathfrak{R}_{0}$ & $I$ &  \\\hline
  $\mathfrak{R}_{2}$ & $D_1,D_2, D_3$ &
  \\\hline
  $\mathfrak{R}_{4}$ & $D_1^2, D_2^2, D_3^2, D_1D_2$ & $E, F$
  \\\hline
  $\mathfrak{R}_{6}$ & $D_1^3, D_2^3, D_3^3, D_1^2D_2$ & $D_1E, D_2F$
  \\\hline
  $\mathfrak{R}_{8}$ & $D_1^4, D_2^4, D_3^4, D_1^3D_2$ & $D_1^2E, D_2^2F$ \\
  \hline
  $\vdots$ & $\vdots$&$\vdots$\\\hline
   $\mathfrak{R}_{2i}$ & $D_1^i, D_2^i, D_3^i, D_1^{i-1}D_2$ & $D_1^{i-2}E, D_2^{i-2}F$
\\
  \hline
\end{tabular}
\end{center}

This is, of course, not the only possible choice of a basis.

The operators $E$ and $F$ themselves are given by
\begin{equation*}
   E= G_4(t)\frac{d^4}{dt^4}+G_3(t)\frac{d^3}{dt^3}+G_2(t)\frac{d^2}{dt^2}+ G_1(t)\frac{d}{dt}+ G_0(t),
\end{equation*}
\begin{equation*}
   F= H_4(t)\frac{d^4}{dt^4}+H_3(t)\frac{d^3}{dt^3}+H_2(t)\frac{d^2}{dt^2}+ H_1(t)\frac{d}{dt}+ H_0(t),
\end{equation*}
with $G_4(t)$ and $H_4(t)$ given by:

$$G_4(t)=\D\frac{(\beta-k_1+2)(k_1-k_2)}{(\beta-k_2+1)}t^2\left(
  \begin{array}{cccc}
0 & 0 & 0 & 0\\
                  -\frac{1}{k_1-k_2-2}(1-t) & 0 &
                  \frac{1}{k_1-k_2-2}(1-t)^2& 0 \\
                 0 &  0&  0&0 \\
                 -\frac{1}{k_1-k_2-1} & 0 &
                 -\frac{1}{k_1-k_2-2}(1-t) & 0\\
  \end{array}
\right)$$

$$H_4(t)=\D\frac{(\beta-k_2+1)(k_1-k_2-2)}{(\beta-k_1+2)}t^2\left(
               \begin{array}{cccc}
                 0 & 0 & 0 & 0\\
                  0 & 0 &0 & 0 \\
                  \frac{1}{k_1-k_2}(1-t)& \frac{1}{k_1-k_2}(1-t)^2 &  0&0 \\
                 \frac{1}{k_2-k_1+1}& \frac{1}{k_2-k_1+1}(1-t) &
                  0& 0\\
               \end{array}
             \right).
$$

The remaining coefficients are too unpleasant to be displayed here.

The operators $D_1, D_2$ and $D_3$ are symmetric with respect to
$W(t)$ at the end of section \ref{sec1}, whereas $E$ and $F$ can be
seen \textbf{not} to be symmetric by considering the coefficients
$G_4(t)$ and $H_4(t)$. The symmetry of $E$ and $F$ would require
$$
WG_4=G_4^*W \quad \mbox{and} \quad WH_4=H_4^*W,
$$
which it is easily seen to be false.

In contrast to the one step example, in this case we have found lots
of relations among the \textbf{generators}
$$
\{I, D_1, D_2, D_3, E, F \}.
$$
We display a few of them below:

$$
D_2\big(D_1+D_3-k_1(\alpha+\beta-k_1+3)\big)=\Theta
$$
\begin{eqnarray*}
\lefteqn{\big[(k_1-k_2)D_2+(k_1-k_2-1)D_3\big]\big[-D_1+(k_1-k_2-1)D_2}\\
\noalign{\smallskip} & &\qquad\qquad
+(k_1-k_2-2)D_3+(1+k_2)(\alpha+\beta-k_2+2)\big]=\Theta
\end{eqnarray*}
$$
D_1 E+E D_3=\big(k_1(\alpha+\beta-k_1+2)+1+k_2\big)E
$$
$$
ED_1-D_1E=(k_1-k_2-1)E
$$
$$
FD_1+D_3F=\big(k_1(\alpha+\beta-k_1+2)+1+k_2\big)F
$$
$$
D_1F-FD_1=(k_1-k_2-1)F
$$
$$
D_2E=\Theta\quad \mbox{and} \quad FD_2=\Theta.
$$

The expressions of $EF$ and $FE$ are obtained in terms of the
elements $\{I, D_2, D_3\}$. For instance, we have:

\begin{eqnarray*}
\lefteqn{\D\frac{(k_1-k_2-1)^3}{(k_1-k_2)^2}FE=D_2\big[(k_1-k_2)D_2+(k_1-k_2-1)(\alpha+\beta-2k_1+3)\big]}\\
\noalign{\smallskip} & &
\big[(k_2-k_1)D_2-(k_1-k_2-1)^2(\alpha+\beta-k_1-k_2+1)\big]\big[D_2+D_3-(\alpha+\beta-k_1-k_2+2)\big]
\end{eqnarray*}

The set of all relations among the generators (something we have not
found) should play the role of the Burchnall--Chaundy curve.

We close our look at this two steps example by stating the
\textbf{conjecture} that $\mathcal{D}$ coincides with the subalgebra
$\mathcal{A}$ generated by $\{I, D_1, D_2, D_3, E, F \}$.

\appendix
\section{Appendix}

\subsection{Explicit expression of $\Pn$ for the one step example}\label{app1}

The polynomials introduced in subsection \ref{sec22} satisfy
$D_1(P_n^*)=P_n^*\Lambda_n(D_1)$  and can be expressed, after an
appropriate choice of the 3 parameters, as
$$
P_n^*(t)=\sum_{l=0}^n A_{l}^nt^{n-l}
$$
where
$$
A_{l}^n=\left(
          \begin{array}{ccc}
            x_{11} & x_{12} & x_{13} \\
            x_{21} & x_{22} & x_{23} \\
            x_{31} & x_{32} & x_{33} \\
          \end{array}
        \right)
$$
and $x_{ij}$ given by:
\begin{eqnarray*}
\lefteqn{x_{11}=(-1)^l\Gamma(n+1)\Gamma(l-2-\alpha-n)\Gamma(-2-\beta-\alpha-2n)(5n^2+2n+2k-2l+2ln^2+2ln}\\
\noalign{\smallskip} &&\quad
-3kl^2+5kl+k^2l^2-3k^2l+7nk+2k^2+2l^2-l\beta^2-3l\beta+3l^2\beta+4kl\beta+3k^2n+n\alpha^2\\
\noalign{\smallskip} & &\quad
+3n\alpha+5n^2\alpha+k^2n^2+7kn^2+2kn^3+n^2\alpha^2+2n^3\alpha-2kln^2-2k^2nl+4\alpha
kn^2+2\alpha
k^2n\\
\noalign{\smallskip} & &\quad +2\alpha^2kn+2ln^2\beta-2\alpha
knl+2\alpha ln\beta+2knl\beta+n^4+2nl\beta+8\alpha kn+2\alpha
ln+4n^3+2\alpha kl\\
\noalign{\smallskip} & &\quad +3\alpha k^2+3\alpha
k+\alpha^2k^2-2\alpha k^2l+\alpha^2k+l^2\beta^2-2kl^2\beta+2\alpha
kl\beta)/\\
\noalign{\smallskip} & &\quad
(\Gamma(l+1)\Gamma(n-l+1)\Gamma(-\alpha-n)\Gamma(l-2-\beta-\alpha-2n)(k+n+1)(k+n))
\end{eqnarray*}
$$
x_{12}=\mbox{{\footnotesize $-\D\frac{(-1)^l(n\alpha+\alpha
k+l\beta+2l+nk-kl+2k+n^2+2n)\Gamma(-\beta-3-\alpha-2n)\Gamma(l-2-\alpha-n)\Gamma(n+1)}{(k+n)\Gamma(n-l+1)\Gamma(l-3-\beta-\alpha-2n)\Gamma(-\alpha-n-1)\Gamma(l+1)}$}}
$$
$$
x_{13}=\frac{(-1)^l\Gamma(n+1)\Gamma(l-2-\alpha-n)\Gamma(-\beta-4-\alpha-2n)}{\Gamma(l+1)\Gamma(n-l+1)\Gamma(-\alpha-n-2)\Gamma(l-4-\beta-\alpha-2n)}
$$
$$
x_{21}=\mbox{{\footnotesize $\D\frac{(-1)^l(n\alpha+\alpha
k+l\beta+2l-\beta-2+nk-kl+2k+n^2+n)\Gamma(-2-\beta-\alpha-2n)\Gamma(l-2-\alpha-n)\Gamma(n+1)}{(k+n)(k+n+1)\Gamma(n-l+1)\Gamma(l-3-\beta-\alpha-2n)\Gamma(-\alpha-n)\Gamma(l)}$}}
$$
\begin{eqnarray*}
\lefteqn{x_{22}=(-1)^{l}\Gamma(n+1)\Gamma(l-2-\alpha-n)\Gamma(-\beta-3-\alpha-2n)(-4lnk+4n^2+4n+4k+12l-ln^2}\\
\noalign{\smallskip} &&\quad
+4ln+2kl^2-10kl+k^2l+2nk-2k^2-4l^2+2k\beta+2l\beta^2+10l\beta-2l^2\beta+2l\alpha\beta+\alpha
k\beta\\
\noalign{\smallskip} &&\quad
-3kl\beta+n^2\beta-k^2n+n\alpha^2+2n\beta+4n\alpha+2n^2\alpha+3nl\beta+kn\beta+\alpha
kn-\alpha ln+n\alpha\beta\\
\noalign{\smallskip} &&\quad +n^3-3\alpha kl-\alpha k^2+4\alpha
k+4l\alpha+\alpha^2k)/\\
\noalign{\smallskip} &&\quad
(2\Gamma(l+1)\Gamma(n-l+1)\Gamma(-\alpha-n-1)\Gamma(l-3-\beta-\alpha-2n)(\beta+1-k)(k+n))
\end{eqnarray*}
$$
x_{23}=-\frac{(-1)^l(\beta+3+\alpha-k+n-l)\Gamma(-\beta-4-\alpha-2n)\Gamma(l-2-\alpha-n)\Gamma(n+1)}{(\beta+1-k)\Gamma(l-4-\beta-\alpha-2n)\Gamma(-\alpha-n-2)\Gamma(n-l+1)\Gamma(l+1)}
$$
$$
x_{31}=\frac{(-1)^l\Gamma(n+1)\Gamma(l-2-\alpha-n)\Gamma(-2-\beta-\alpha-2n)}{(k+n)(k+n+1)\Gamma(l-1)\Gamma(n-l+1)\Gamma(-\alpha-n)\Gamma(l-4-\beta-\alpha-2n)}
$$
$$
x_{32}=\frac{(-1)^l(\beta+3+\alpha-k+n-l)\Gamma(-\beta-3-\alpha-2n)\Gamma(l-2-\alpha-n)\Gamma(n+1)}{(k+n)(\beta+1-k)\Gamma(l-4-\beta-\alpha-2n)\Gamma(n-l+1)\Gamma(-\alpha-n-1)\Gamma(l)}
$$
$$
x_{33}=\mbox{{\footnotesize
$\D\frac{(-1)^l(\beta+4+\alpha-k+n-l)(\beta+3+\alpha-k+n-l)\Gamma(-\beta-4-\alpha-2n)\Gamma(l-2-\alpha-n)\Gamma(n+1)}{(\beta+1-k)(\beta+2-k)\Gamma(l-4-\beta-\alpha-2n)\Gamma(-\alpha-n-2)\Gamma(n-l+1)\Gamma(l+1)}$}}
$$
Here, $\Gamma$ denotes the Gamma Euler function.

\subsection{Some constants}\label{app2}

Here we write down the values of $C_{22}, C_{33}, C_{44}$ of
subsection \ref{sec31}:

\begin{eqnarray*}
\lefteqn{C_{22}=(18+k_1^2\alpha\beta-k_2^2\alpha^2+21\alpha-k_2^2\alpha\beta+9k_2\alpha\beta-k_1^2k_2\alpha+k_2\alpha^3+k_1k_2^2\alpha+k_2\alpha\beta^2-k_1\alpha\beta^2}\\
\noalign{\smallskip} & &\quad
+2k_2\alpha^2\beta-2k_1\alpha^2\beta-k_1k_2\alpha-11k_1\alpha\beta+k_1^2\alpha^2-27k_1+12\beta+15k_2-5k_2^2+8\alpha^2+8k_1^2\\
\noalign{\smallskip} & &\quad
+\alpha^3-k_1\alpha^3-k_1k_2-27k_1\alpha+6k_1^2\alpha-4k_2^2\alpha-9k_1\alpha^2+7k_2\alpha^2-2k_1^2k_2+2k_1k_2^2+17k_2\alpha\\
\noalign{\smallskip} & &\quad
+2\beta^2+\alpha\beta^2+10\alpha\beta+2\alpha^2\beta-2k_1\beta^2-15k_1\beta+2k_2\beta^2+11k_2\beta+2k_1^2\beta-2k_2^2\beta)/\\
\noalign{\smallskip} & &\quad
((k_2-k_1+1)(\alpha+\beta-k_1+3)(\alpha+\beta-k_2+3))
\end{eqnarray*}

\begin{eqnarray*}
\lefteqn{C_{33}=(16+k_1^2\alpha\beta-k_2^2\alpha^2+18\alpha-k_2^2\alpha\beta+9k_2\alpha\beta-k_1^2k_2\alpha+k_2\alpha^3+k_1k_2^2\alpha+k_2\alpha\beta^2-k_1\alpha\beta^2}\\
\noalign{\smallskip} & &\quad
+2k_2\alpha^2\beta-2k_1\alpha^2\beta+5k_1k_2\alpha-11k_1\alpha\beta+k_1^2\alpha^2-18k_1+14\beta+4k_2-10k_2^2+8\alpha^2+3k_1^2\\
\noalign{\smallskip} & &\quad
+\alpha^3-k_1\alpha^3+9k_1k_2-21k_1\alpha+3k_1^2\alpha-7k_2^2\alpha-9k_1\alpha^2+7k_2\alpha^2-2k_1^2k_2+2k_1k_2^2+11k_2\alpha\\
\noalign{\smallskip} & &\quad
+2\beta^2+\alpha\beta^2+10\alpha\beta+2\alpha^2\beta-2k_1\beta^2-15k_1\beta+2k_2\beta^2+11k_2\beta+2k_1^2\beta-2k_2^2\beta)/\\
\noalign{\smallskip} & &\quad
((k_2-k_1+1)(\alpha+\beta-k_1+4)(\alpha+\beta-k_2+2))
\end{eqnarray*}

\begin{eqnarray*}
\lefteqn{C_{44}=(24+23\alpha-k_2\alpha\beta+k_1k_2\alpha-k_1\alpha\beta-7k_1+17\beta-10k_2+8\alpha^2+\alpha^3+3k_1k_2}\\
\noalign{\smallskip} & &\quad
-5k_1\alpha-k_1\alpha^2-k_2\alpha^2-6k_2\alpha+3\beta^2+\alpha\beta^2+11\alpha\beta+2\alpha^2\beta-3k_1\beta-3k_2\beta)\\
\noalign{\smallskip} & &\quad
/((\alpha+\beta-k_1+4)(\alpha+\beta-k_2+2))
\end{eqnarray*}

\subsection{Explicit expression of $\Pn$ for the two steps example}\label{app3}

The polynomials introduced in subsection \ref{sec32} satisfy
$D_1(P_n^*)=P_n^*\Lambda_n(D_1)$  and can be expressed, after an
appropriate choice of the 4 parameters, as
$$
P_n^*(t)=\sum_{l=0}^n A_{l}^nt^{n-l}
$$
where
$$
A_{l}^n=\left(
          \begin{array}{cccc}
            x_{11} & x_{12} & x_{13} &  x_{14}\\
            x_{21} & x_{22} & x_{23} &  x_{24}\\
            x_{31} & x_{32} & x_{33} &  x_{34}\\
            x_{41} & x_{42} & x_{43} &  x_{44}\\
          \end{array}
        \right)
$$
and $x_{ij}$ given by:
\begin{eqnarray*}
\lefteqn{x_{11}=(-1)^l\Gamma(n+1)\Gamma(l-2-\alpha-n)\Gamma(-2-\beta-\alpha-2n)(-ln^2k_1-ln^2k_2+\alpha
k_2l\beta-k_2l^2\beta}\\
\noalign{\smallskip} &&\quad -k_1l^2\beta+5n^2+2k_1l\beta+\alpha
k_1l\beta+2k_2l\beta+2nk_2+5nk_1+2k_1k_2-l\beta^2+3l^2\beta+3n^2k_2\\
\noalign{\smallskip} & &\quad
+4n^2k_1-3l\beta+3nk_1k_2-k_1l^2-2l+2n+3\alpha
k_1+2k_1+2l^2-2k_2l^2+4k_2l+k_1l+lnk_2\\
\noalign{\smallskip} & &\quad
-lnk_1+2ln+2ln^2+k_1k_2l^2-3k_1k_2l+l^2\beta^2+\alpha^2k_1+\alpha^2k_1k_2+3\alpha
k_1k_2+k_1n^3+2\alpha k_2l\\
\noalign{\smallskip} & &\quad +2\alpha k_1k_2n-\alpha k_1nl-\alpha
k_2nl+k_1nl\beta+k_2n^3+5n^2\alpha+2nl\beta+4n^3+k_1k_2n^2+k_2nl\beta\\
\noalign{\smallskip} & &\quad +3\alpha k_2n+2\alpha k_2n^2+5\alpha
k_1n+2\alpha k_1n^2+2\alpha
nl+\alpha^2k_1n+n\alpha^2k_2+3n\alpha+n\alpha^2-2k_1k_2nl\\
\noalign{\smallskip} & &\quad -2\alpha
k_1k_2l+2\alpha nl\beta+2n^3\alpha+n^4+2ln^2\beta+n^2\alpha^2)/\\
\noalign{\smallskip} & &\quad
(\Gamma(l+1)\Gamma(n-l+1)\Gamma(-\alpha-n)\Gamma(l-2-\beta-\alpha-2n)(n+k_2+1)(k_1+n))
\end{eqnarray*}
$$
x_{12}=\mbox{{\footnotesize
$\D\frac{(-1)^{l+1}\Gamma(n+1)\Gamma(l-2-\alpha-n)\Gamma(-\beta-3-\alpha-2n)(2l+2k_1+\alpha
k_1+l\beta+n^2+n\alpha+2n-k_1l+nk_1)}{\Gamma(l+1)\Gamma(n-l+1)\Gamma(-\alpha-n-1)\Gamma(l-3-\beta-\alpha-2n)(k_1+n)}$}}
$$
$$
x_{13}=\mbox{{\scriptsize
$\D\frac{(-1)^{l+1}\Gamma(n+1)\Gamma(l-2-\alpha-n)\Gamma(-\beta-3-\alpha-2n)(2k_2+2+\alpha+l+n^2+3n+\alpha
k_2+nk_2+n\alpha+l\beta-k_2l)}{(n+k_2+1)\Gamma(l-3-\beta-\alpha-2n)\Gamma(-\alpha-n-1)\Gamma(n-l+1)\Gamma(l+1)}$}}
$$
$$
x_{14}=\frac{(-1)^l\Gamma(n+1)\Gamma(l-2-\alpha-n)\Gamma(-\beta-4-\alpha-2n)}{\Gamma(l+1)\Gamma(n-l+1)\Gamma(-\alpha-n-2)\Gamma(l-4-\beta-\alpha-2n)}
$$
$$
x_{21}=\mbox{{\scriptsize
$\D\frac{(-1)^l\Gamma(n+1)\Gamma(l-2-\alpha-n)\Gamma(-2-\beta-\alpha-2n)(n\alpha+\alpha
k_1+l\beta+2l-\beta-2+nk_1-k_1l+2k_1+n^2-n)}{(k_1+n)(n+k_2+1)\Gamma(l)\Gamma(n-l+1)\Gamma(-\alpha-n)\Gamma(l-3-\beta-\alpha-2n)}$}}
$$
\begin{eqnarray*}
\lefteqn{x_{22}=(-1)^l\Gamma(n+1)\Gamma(l-2-\alpha-n)\Gamma(-\beta-3-\alpha-2n)(ln^2k_1-ln^2k_2+k_2^2k_1l+2nk_1^2l-2k_2^2l}\\
\noalign{\smallskip} &&\quad -k_1^2k_2l+\alpha
k_2l\beta-k_2l^2\beta+k_1l^2\beta+4n^2+\alpha
k_1\beta+2k_1k_2\beta-k_2^2l\beta-k_1l\beta^2-6k_1l\beta-\alpha k_1l\beta\\
\noalign{\smallskip} &&\quad +\alpha
k_1k_2\beta+k_1^2l\beta+k_2l\beta^2+3k_2l\beta+2nk_2+2k_1k_2-4k_1^2-k_1^2\beta
n+2l\beta^2+2\alpha l\beta-2l^2\beta\\
\noalign{\smallskip} &&\quad
+3n^2k_2-2nk_2^2-3n^2k_1+10l\beta+2k_1\beta-k_1^2l^2+4k_1^2l-2k_2^2k_1-4nk_1^2+2k_1^2k_2+5nk_1k_2\\
\noalign{\smallskip} &&\quad +4k_1l^2+12l+4n+4\alpha
k_1+4\alpha l+4k_1-4l^2-2k_1^2\beta-2k_2l^2+2k_2l-12k_1l-4lnk_1\\
\noalign{\smallskip} &&\quad
+4ln-ln^2+k_1k_2l^2-k_1k_2l-\alpha^2k_1^2-k_1^2\beta\alpha+\alpha^2k_1-4\alpha
k_1^2+\alpha^2k_1k_2+3\alpha k_1k_2+\alpha k_1^2k_2\\
\noalign{\smallskip} &&\quad -\alpha k_2^2k_1-k_1n^3-5\alpha
k_1l+2\alpha k_2l+k_1k_2n\beta+3\alpha k_1k_2n+\alpha k_1nl-\alpha
k_2nl+n^2\beta-k_1nl\beta\\
\noalign{\smallskip} &&\quad -\alpha k_1n\beta+\alpha
k_2n\beta+k_2n^3+2n^2\alpha-k_1^2n^2-k_2^2n^2+2n\beta-k_2^2k_1n+k_1^2k_2n+3nl\beta-k_1n^2\beta\\
\noalign{\smallskip} &&\quad
-k_1n\beta+k_2n^2\beta+2k_2n\beta+n^3+2k_1k_2n^2+k_2nl\beta+3\alpha
k_2n+2\alpha k_2n^2-2\alpha k_1n-2\alpha k_1n^2\\
\noalign{\smallskip} &&\quad -n\alpha k_2^2-2\alpha k_1^2n-\alpha
nl-\alpha^2k_1n+n\alpha^2k_2+4n\alpha+n\alpha^2+n\alpha\beta-2k_1k_2nl+2\alpha
k_1^2l-2\alpha
k_1k_2l)/\\
\noalign{\smallskip} &&\quad
(\Gamma(l+1)\Gamma(n-l+1)\Gamma(-\alpha-n-1)\Gamma(l-3-\beta-\alpha-2n)(\beta+1-k_2)(k_1+n)(k_2-k_1+2))
\end{eqnarray*}
$$
x_{23}=\frac{(-1)^{l+1}\Gamma(n+1)\Gamma(l-2-\alpha-n)\Gamma(-\beta-3-\alpha-2n)}{(n+k_2+1)\Gamma(l)\Gamma(n-l+1)\Gamma(-\alpha-n-1)\Gamma(l-4-\beta-\alpha-2n)}
$$
$$
x_{24}=\frac{(-1)^{l+1}\Gamma(n+1)\Gamma(l-2-\alpha-n)\Gamma(-\beta-4-\alpha-2n)(\beta+3+\alpha-k_2+n-l)}{(\beta+1-k_2)\Gamma(l-4-\beta-\alpha-2n)\Gamma(-\alpha-n-2)\Gamma(n-l+1)\Gamma(l+1)}
$$
$$
x_{31}=\mbox{{\footnotesize
$\D\frac{(-1)^l\Gamma(n+1)\Gamma(l-2-\alpha-n)\Gamma(-2-\beta-\alpha-2n)(n\alpha+\alpha+\alpha
k_2+l\beta+l-\beta+nk_2-k_2l+2k_2+n^2+2n)}{(k_1+n)(n+k_2+1)\Gamma(l)\Gamma(n-l+1)\Gamma(-\alpha-n)\Gamma(l-3-\beta-\alpha-2n)}$}}
$$
$$
x_{32}=\frac{(-1)^{l+1}\Gamma(n+1)\Gamma(l-2-\alpha-n)\Gamma(-\beta-3-\alpha-2n)}{(k_1+n)\Gamma(l)\Gamma(n-l+1)\Gamma(-\alpha-n-1)\Gamma(l-4-\beta-\alpha-2n)}
$$
\begin{eqnarray*}
\lefteqn{x_{33}=(-1)^{l+1}\Gamma(n+1)\Gamma(l-2-\alpha-n)\Gamma(-\beta-3-\alpha-2n)(6+ln^2k_1-2lnk_2^2-ln^2k_2+k_2^2k_1l}\\
\noalign{\smallskip} &&\quad -5k_2^2l-k_1^2k_2l+\alpha
k_2l\beta-k_2l^2\beta+k_1l^2\beta+6n^2-\alpha k_1\beta+\alpha
k_2^2\beta-2k_1k_2\beta+2\alpha
k_2\beta-k_2^2l\beta\\
\noalign{\smallskip} &&\quad -k_1l\beta^2-5k_1l\beta-\alpha
k_1l\beta-\alpha
k_1k_2\beta+k_1^2l\beta+k_2l\beta^2+4k_2l\beta+\alpha^2+16nk_2-14nk_1-10k_1k_2\\
\noalign{\smallskip} &&\quad
+2k_1^2+6k_2^2+\alpha\beta+7n^2k_2+5nk_2^2-7n^2k_1-l\beta-2k_1\beta+k_1^2l-2k_2^2k_1+3nk_1^2+2k_1^2k_2\\
\noalign{\smallskip} &&\quad
+k_2^2l^2+2k_2^2\beta-9nk_1k_2+k_1l^2-3l+4k_2\beta+11n+10\alpha
k_2-6\alpha
k_1-\alpha l+12k_2-8k_1+5\alpha\\
\noalign{\smallskip} &&\quad
-k_2l^2-2k_1l-2lnk_2+2lnk_1-4ln-ln^2-k_1k_2l^2+5k_1k_2l+\alpha^2k_2^2-\alpha^2k_1+2\alpha^2k_2+\alpha
k_1^2\\
\noalign{\smallskip} &&\quad +5\alpha k_2^2-\alpha^2k_1k_2-7\alpha
k_1k_2+\alpha k_1^2k_2+2\beta-\alpha
k_2^2k_1-2\alpha k_2^2l-k_1n^3-\alpha k_2l-k_1k_2n\beta\\
\noalign{\smallskip} &&\quad -3\alpha k_1k_2n+\alpha k_1nl-\alpha
k_2nl+n^2\beta-k_1nl\beta-\alpha
k_1n\beta+\alpha k_2n\beta+k_2n^3+2n^2\alpha+k_1^2n^2+k_2^2n^2\\
\noalign{\smallskip} &&\quad
+3n\beta-k_2^2k_1n+k_1^2k_2n+k_2^2n\beta-nl\beta-k_1n^2\beta-3k_1n\beta+k_2n^2\beta+4k_2n\beta+n^3-2k_1k_2n^2\\
\noalign{\smallskip} &&\quad +k_2nl\beta+9\alpha k_2n+2\alpha
k_2n^2-8\alpha k_1n-2\alpha k_1n^2+2n\alpha
k_2^2+\alpha k_1^2n-\alpha nl-\alpha^2k_1n+n\alpha^2k_2\\
\noalign{\smallskip} &&\quad
+7n\alpha+n\alpha^2+n\alpha\beta+2k_1k_2nl+2\alpha k_1k_2l)/\\
\noalign{\smallskip} &&\quad
(\Gamma(l+1)\Gamma(n-l+1)\Gamma(-\alpha-n-1)\Gamma(l-3-\beta-\alpha-2n)(n+k_2+1)(k_1-k_2)(\beta+2-k_1))\\
\end{eqnarray*}
$$
x_{34}=\frac{(-1)^{l+1}\Gamma(n+1)\Gamma(l-2-\alpha-n)\Gamma(-\beta-4-\alpha-2n)(\beta+4+\alpha-k_1+n-l)}{\Gamma(l+1)\Gamma(n-l+1)\Gamma(-\alpha-n-2)\Gamma(l-4-\beta-\alpha-2n)(\beta+2-k_1)}
$$
$$
x_{41}=\frac{(-1)^l\Gamma(n+1)\Gamma(l-2-\alpha-n)\Gamma(-2-\beta-\alpha-2n)}{(k_1+n)(n+k_2+1)\Gamma(l-1)\Gamma(n-l+1)\Gamma(-\alpha-n)\Gamma(l-4-\beta-\alpha-2n)}
$$
$$
x_{42}=\frac{(-1)^l(\beta+3+\alpha-k_2+n-l)\Gamma(n+1)\Gamma(l-2-\alpha-n)\Gamma(-\beta-3-\alpha-2n)}{(k_1+n)(\beta+1-k_2)\Gamma(l)\Gamma(n-l+1)\Gamma(-\alpha-n-1)\Gamma(l-4-\beta-\alpha-2n)}
$$
$$
x_{43}=\frac{(-1)^l(\beta+4+\alpha-k_1+n-l)\Gamma(n+1)\Gamma(l-2-\alpha-n)\Gamma(-\beta-3-\alpha-2n)}{(n+k_2+1)(\beta+2-k_1)\Gamma(l)\Gamma(n-l+1)\Gamma(-\alpha-n-1)\Gamma(l-4-\beta-\alpha-2n)}
$$
$$
x_{44}=\mbox{{\footnotesize
$\D\frac{(-1)^l\Gamma(n+1)\Gamma(l-2-\alpha-n)\Gamma(-\beta-4-\alpha-2n)(\beta+3+\alpha-k_2+n-l)(\beta+4+\alpha-k_1+n-l)}{\Gamma(l+1)\Gamma(n-l+1)\Gamma(-\alpha-n-2)\Gamma(l-4-\beta-\alpha-2n)(\beta+1-k_2)(\beta+2-k_1)}$}}
$$

\end{document}